\documentclass[11pt]{article}
\usepackage{amsmath,amssymb,amsthm}
\usepackage{mathtools}
\usepackage{tikz-cd} 
\usepackage{hyperref} 
\theoremstyle{definition} 
\newtheorem{definition}{Definition}[section]
\newtheorem{remark}[definition]{Remark}
\theoremstyle{plain}
\newtheorem{theorem}[definition]{Theorem}
\newtheorem{corollary}[definition]{Corollary}
\newtheorem{proposition}[definition]{Proposition}
\newtheorem{lemma}[definition]{Lemma}
\newtheorem{example}[definition]{Example} \title{Bott--Chern and Aeppli Cohomology of toroidal groups} 
\author{Xiaojun Wu} 
\date{\today} 
\begin{document} 
\def\A{\mathcal{A}}
\def\cI{\mathcal{I}}
\def\Z{\mathbb{Z}}
\def\Q{\mathbb{Q}}  \def\C{\mathbb{C}}
 \def\R{\mathbb{R}}
 \def\N{\mathbb{N}}
 \def\H{\mathbb{H}}
  \def\P{\mathbb{P}}
  \def\dist{\mathrm{dist}}
  \def\d{\partial}
 \def\dbar{{\overline{\partial}}}
\def\dzbar{{\overline{dz}}}
\def \ddbar {\partial \overline{\partial}}
\def\cB{\mathcal{B}}
\def\cD{\mathcal{D}}  \def\cO{\mathcal{O}}
\def\cbarO{\overline{\mathcal{O}}}
\def\D{\mathcal{D}}
\def\cC{\mathcal{C}}
\def\cF{\mathcal{F}}
\def \rank{\mathrm{rank}}
\def \deg{\mathrm{deg}}
\def \Ker{\mathrm{Ker}}
\def \delbar{\overline{\partial}}
\def \tot{\mathrm{Tot}}
\def \id{\mathrm{id}}
\bibliographystyle{plain}
\def \End{\mathrm{End}}
\def \dim{\mathrm{dim}}
\def \div{\mathrm{div}}
\def \ker{\mathrm{Ker}}
\def \im{\mathrm{Im}}
\def \rC{\mathrm{Cone}}
\newcommand{\Ub}{\mathcal{U}}
\newcommand{\dcech}{\check{\delta}}
\newcommand{\dcechg}{\delta\!\!\!\check{\delta}}
\newcommand{\lc}{\mathcal{L}}
\newcommand{\ec}{\mathcal{E}}
\maketitle 
\begin{abstract}
We compute the Bott--Chern and Aeppli cohomology groups of theta toroidal groups in arbitrary dimensions, obtaining explicit dimension formulas and bases in all bidegrees. For general toroidal groups, we determine the Hausdorff reductions of these cohomology groups and show that they coincide with those of theta toroidal groups. We also compute the Bott--Chern and Aeppli cohomology groups themselves in certain specific bidegrees, without passing to their Hausdorff reductions.
 \end{abstract}
\section{Introduction} Bott--Chern and Aeppli cohomology provide refinements of Dolbeault cohomology that  capture finer analytic and non-Kähler features of complex manifolds. Since their introduction by Bott--Chern \cite{BC65} and Aeppli \cite{Aep62}, these cohomology theories have played an important role in complex geometry, particularly in the study of non-compact or non-K\"ahler complex manifolds. 

A sheaf-theoretic interpretation due to Demailly \cite{agbook} and further developed by Schweitzer \cite{Sch07} realizes Bott--Chern and Aeppli cohomology groups as hypercohomology groups of certain complexes, thereby providing a natural framework for computations once the de Rham and Dolbeault cohomologies are known. 

Toroidal groups form a distinguished class of connected complex Lie groups. By definition, a toroidal group admits no nonconstant holomorphic functions. Equivalently, it can be realized as a quotient \( X = \mathbb{C}^n / \Gamma, \) where $\Gamma \subset \mathbb{C}^n$ is a discrete subgroup of real rank $n+m$ with $0 \le m < n$, satisfying a Diophantine-type irrationality condition \cite{Kop,AK}. The geometry and cohomology of toroidal groups exhibit features that sharply distinguish them from compact complex manifolds. In particular, the Dolbeault cohomology groups of toroidal groups need not be Hausdorff in their natural topologies, and their dimensions may be infinite in general. A fundamental breakthrough due to Kazama \cite{Kazama} shows that under a Diophantine condition on the lattice $\Gamma$, the resulting toroidal group belongs to the class of \emph{theta toroidal groups}. For such groups, the Dolbeault cohomology groups are finite-dimensional and generated by invariant constant-coefficient forms.
Moreover, de Rham cohomology is also well-known \cite{Ume}. 
The remaining toroidal groups are called \emph{wild toroidal groups}, i.e.\ complex Lie groups of the form $X = \mathbb{C}^n / \Gamma$, where $\Gamma \subset \mathbb{C}^n$ is a discrete subgroup of real rank $n + m$ with $0 \le m < n$, and for which some Dolbeault cohomology groups are infinite-dimensional.
The Dolbeault cohomology of such toroidal groups is described in \cite{FRR19}.

Unlike the Dolbeault cohomology, the Bott--Chern and Aeppli cohomology of theta toroidal groups are more complicated and subtle in their structure. 
It is well-known that, as a smooth manifold, a theta toroidal group is diffeomorphic to a product of circles and Euclidean lines. This diffeomorphism allows one to consider formal Fourier expansions of smooth forms on the toroidal group. 
Using such expansions, one can explicitly construct a family of differential forms whose coefficients are constant with respect to the global toroidal coordinates. These forms are invariant under the natural toroidal group action and hence can be regarded as smooth forms on the toroidal group. 
When one examines the Dolbeault cohomology, these invariant constant-coefficient forms provide  natural candidates for generators. To verify that these forms represent linearly independent cohomology classes, one can use dimension-counting arguments to show that the Hodge--de~Rham spectral sequence degenerates at \(E_1\) and that these forms provide a basis of the Dolbeault cohomology groups. 
In contrast, for Bott--Chern or Aeppli cohomology, there is no known natural spectral sequence analogous to the Hodge--de~Rham spectral sequence. Consequently, it is not straightforward to determine whether the cohomology classes represented by invariant constant-coefficient forms form a basis of the Bott--Chern or Aeppli cohomology groups. In fact, the Aeppli cohomology groups may admit additional generators arising from non-constant pluriharmonic functions.
\paragraph{}
The first aim of this paper is to compute explicitly the Bott--Chern and Aeppli cohomology groups of theta toroidal groups, using their sheaf-theoretic descriptions for all complex dimensions. 

We first carry out a complete computation in complex dimension two. This recovers and rederives results of Tanaka \cite{Tan25} by a method based on hypercohomology and mapping cone constructions. Tanaka’s original approach in \cite{Tan25} relies on analytic techniques. More precisely, a two-dimensional theta toroidal group is diffeomorphic to \((S^1)^3 \times \R\). After performing Fourier expansion with respect to the \((S^1)^3\)-variables, the equations associated with the operators \(\d\), \(\dbar\), and \(\d\dbar\) acting on each Fourier coefficient reduce to ordinary differential equations. These equations can be solved explicitly, which allows a direct computation of the relevant cohomology groups.
However, extending this method to higher dimensions seems difficult. Indeed, the resulting equations would form a system of linear partial differential equations with constant coefficients, rather than a system of linear ordinary differential equations with constant coefficients. In the former case, there is no general theory that readily provides the dimension of the solution space.
To illustrate how our cohomological approach can circumvent this analytic difficulty, we make the following observation. In the proof of Theorem \ref{BC-A-theta}, we establish the injectivity of the natural map from Bott--Chern cohomology to de Rham cohomology. This injectivity, in turn, implies the solvability of a certain \(\d\dbar\)-equation.


We then extend these computations to arbitrary-dimensional theta toroidal groups (Theorem \ref{BC-A-theta}). By exploiting the explicit description of de Rham and Dolbeault cohomology in terms of invariant constant-coefficient  forms due to Umeno \cite{Ume} and Kazama \cite{Kazama}, we obtain closed formulas for the dimensions of Bott--Chern and Aeppli cohomology groups in all bidegrees and give explicit bases for these groups. In particular, we prove that all Bott--Chern and Aeppli cohomology groups of theta toroidal groups are finite-dimensional, despite the non-compactness of the underlying manifolds. This work provides a systematic and uniform approach to Bott--Chern and Aeppli cohomology on toroidal groups, clarifying the role of invariant forms and hypercohomology techniques in a setting where standard Hodge-theoretic tools are not available.
\paragraph{} 
The second aim of this paper is to study the Bott--Chern and Aeppli cohomology of wild toroidal groups. For these groups, the Dolbeault, Bott--Chern, and Aeppli cohomology groups naturally arise as quotients of Fréchet spaces. Kazama \cite{Kazama} showed that, equipped with these natural quotient topologies, the Dolbeault cohomology groups in the infinite-dimensional case are non-Hausdorff.

To extract meaningful invariants, we systematically consider the \emph{Hausdorff reduction} of these cohomology groups, obtained by quotienting by the closure of zero. This construction is functorial but not exact, and therefore the hypercohomology  computations for theta toroidal groups do not directly extend to this setting. While the Hausdorff reductions of Dolbeault cohomology are known to be isomorphic to the corresponding theta toroidal groups \cite{FRR19}, comparatively little is known about the Bott--Chern and Aeppli cohomology groups.

In the case of wild toroidal groups, a direct hypercohomology computation cannot be carried out, since the Hausdorff reduction is a non-exact functor when the morphisms are not closed. However, our previous hypercohomology computations of the Bott--Chern and Aeppli cohomology groups for theta toroidal groups provide explicit bases of the corresponding groups (see Remark \ref{Abase2} and Remark \ref{BCbase2}), and the same expressions remain meaningful in the wild case.

Moreover, the proof of \cite{FRR19} provides an explicit method to formally solve the $\partial$- or $\bar{\partial}$-equations on a wild toroidal group. Using a zig-zag argument inspired by the usual zig-zag procedure in spectral sequence computations, we first show that the forms appearing in these explicit bases generate the Hausdorff reduction of the Bott--Chern cohomology. Motivated by the theta-toroidal case, we then prove that the naturally induced map from the Hausdorff reduction of the Bott--Chern cohomology to the de Rham cohomology is injective, thereby establishing the linear independence of the corresponding Hausdorff reduction of the cohomology classes.

In the Aeppli case, we rely on formally solving the equations defining Bott--Chern cohomology. In particular, we conclude in Theorems \ref{BCwild2} and \ref{Awild} that the Hausdorff reductions of the Bott--Chern and Aeppli cohomology groups coincide with those of theta toroidal groups.

Finally, in the last part of this section, we use Malgrange’s vanishing theorem for coherent analytic sheaves on toroidal groups to compute certain cohomology groups in specific bidegrees. A key feature of our approach is that, whenever possible, the statements are established at the level of the cohomology groups themselves, rather than only for their Hausdorff reductions.
 \paragraph{}
The organization of the paper is as follows. In Section 2, we compute the Bott--Chern and Aeppli cohomology groups of a two-dimensional theta toroidal group, thereby rederiving the results of Tanaka. There are two reasons for presenting this explicit calculation. First, Tanaka’s computation determines only the Hausdorff reductions of these cohomology groups. We show that, in this case, the natural maps to the Hausdorff reductions are in fact isomorphisms (Remark \ref{rem-hau-red}). Second, this example already illustrates how to compute the kernels and cokernels of the morphisms appearing in the exact sequences.
In Section 3, using a hypercohomological description, we determine these cohomology groups in all bidegrees for theta toroidal groups in arbitrary dimension. We collect some basic results on toroidal groups in the opening subsections of Sections 3 and 4. In the final section, we discuss the Hausdorff reduction of the cohomology groups associated with general toroidal groups.
 \paragraph{}
\textbf{Acknowledgements.} I would like to thank my postdoctoral advisor, Professor Takayuki Koike, and Dr.~Jinichiro Tanaka for several helpful and stimulating discussions related to this work. This research was supported by the JSPS Postdoctoral Fellowships for Research in Japan (Standard). I am also grateful to Osaka Metropolitan University and the Osaka Central Advanced Mathematical Institute (OCAMI) for providing an excellent research environment.
\section{Two-dimensional case}
\subsection{Preliminaries}
Bott--Chern and Aeppli cohomology are refinements of Dolbeault cohomology; we briefly recall their definitions for completeness (see, for example, \cite{Aep62,BC65,agbook,Sch07}). 
\begin{definition}\label{BC}
 Let $X$ be a complex manifold. The \emph{Bott--Chern cohomology} of $X$ is defined by \[ H^{p,q}_{BC}(X) = \frac{ \ker \partial \cap \ker \bar\partial \subset \mathcal A^{p,q}(X) }{ \operatorname{Im}\,\partial\bar\partial \subset \mathcal A^{p,q}(X) }. \] \end{definition}
 \begin{definition} 
\label{Ae} 
 Let $X$ be a complex manifold. The \emph{Aeppli cohomology} of $X$ is defined by \[ H^{p,q}_{A}(X) = \frac{ \ker \partial\bar\partial \subset \mathcal A^{p,q}(X) }{ \operatorname{Im}\,\partial + \operatorname{Im}\,\bar\partial }. \] \end{definition}
 We now recall the sheaf-theoretic interpretation of Bott--Chern cohomology due to Demailly (see~\cite[Chap.~VI, \S12.1]{agbook}, \cite{Sch07}).
This sheaf-theoretic description provides a natural framework for relating Bott--Chern and Aeppli cohomology to Dolbeault cohomology and de~Rham cohomology.
 \begin{proposition}
\label{Dem-Sch} 
The Bott--Chern cohomology group defined in Definition \ref{BC} is the hypercohomology group \[ H^{p,q}_{BC}(X) = \mathbb H^{p+q}(X,\mathcal B^\bullet_{p,q}) \] of the Bott--Chern complex \[ \mathcal B^\bullet_{p,q} : (2\pi \sqrt{-1})^p\mathbb C \xrightarrow{\Delta} \mathcal O \oplus \bar{\mathcal O} \to \Omega^1 \oplus \bar{\Omega}^1 \to \cdots \to \Omega^{p-1} \oplus \bar{\Omega}^{p-1} \to \cdots \to \bar{\Omega}^{q-1} \to 0, \] 
for $p \leq q$ or
\[ \mathcal B^\bullet_{p,q} : (2\pi \sqrt{-1})^p\mathbb C \xrightarrow{\Delta} \mathcal O \oplus \bar{\mathcal O} \to \Omega^1 \oplus \bar{\Omega}^1 \to \cdots \to \Omega^{q-1} \oplus \bar{\Omega}^{q-1} \to \cdots \to \Omega^{p-1} \to 0, \] 
for $p \geq q$
 where $\Delta$ is multiplication by $1$ on the first component and by $-1$ on the second component.

Moreover, the Aeppli cohomology group defined in Definition \ref{Ae} is the hypercohomology group
\[ H^{p,q}_A(X) = \mathbb H^{p+q+1}(X,\mathcal B^\bullet_{p+1,q+1}). \]  
 \end{proposition} 
The second statement was observed by Schweitzer~\cite[Section 4.d]{Sch07}.  
 
The Bott--Chern complex is quasi-isomorphic to the complex $\mathcal L^\bullet_{p,q}$ introduced in \cite[Chap.~VI, \S12.1]{agbook}, which we recall below.

For $k \le p+q-2$, set \[ \mathcal L^k_{p,q} = \bigoplus_{\substack{r+s=k \\ r<p,\ s<q}} \mathcal E^{r,s}, \] while for $k \ge p+q$, set \[ \mathcal L^{k-1}_{p-1,q-1} = \bigoplus_{\substack{r+s=k \\ r\geq p,\ s\geq q}} \mathcal E^{r,s}. \] The differential of $\mathcal L^\bullet_{p,q}$ is given by \[ \mathcal L^0 \xrightarrow{\;\operatorname{pr}_{\mathcal L^1}\circ d\;} \mathcal L^1 \xrightarrow{\;\operatorname{pr}_{\mathcal L^2}\circ d\;} \cdots \to \mathcal L^{p+q-2} \xrightarrow{\;\frac{\sqrt{-1}}{2\pi}\partial\bar\partial\;} \mathcal L^{p+q-1} \xrightarrow{\;d\;} \mathcal L^{p+q} \xrightarrow{\;d\;} \cdots . \]
Here $\mathcal E^{r,s}$ denotes the sheaf of smooth $(r,s)$-forms on $X$.

\begin{remark}
\label{rem-attention-topo}
Since
the Bott--Chern complex is quasi-isomorphic to the complex $\mathcal L^\bullet_{p,q}$, this quasi-isomorphism induces a linear isomorphism between their hypercohomology groups. However, the hypercohomology of the Bott--Chern complex is not, a priori, equipped with a topological vector space structure. Throughout the paper, whenever we consider the Hausdorff reductions of Bott--Chern and Aeppli cohomology groups, we equip the cohomology groups with the topology induced by the natural topology on the spaces of smooth forms appearing as components of \(\mathcal L^\bullet_{p,q}\).
\end{remark}
\begin{remark}
\label{rem-quasi-isom}
Let $\mathcal A^\bullet_p$ be the complex
\[
\mathcal O_X \to \Omega^1_X \to \cdots \to \Omega^{p-1}_X \to 0,
\]
starting in degree $0$, and let
$\iota_p:\mathcal A^\bullet_p\to\mathcal L^\bullet_{p,q}$
be the natural inclusion. Consider the morphism of complexes
\[
\iota_p+\overline{\iota_q}:
\mathcal A^\bullet_p\oplus\overline{\mathcal A^\bullet_q}
\longrightarrow \mathcal L^\bullet_{p,q}.
\]
Let
$\Delta:\mathbb C\to
\mathcal A^\bullet_p\oplus\overline{\mathcal A^\bullet_q}$
be the morphism of complexes induced by the diagonal map
$c\mapsto(c,-c)$ in degree $0$, where $\mathbb C$ is regarded as a
complex concentrated in degree $0$. We assume that $p,q\geq 1$.
The Bott--Chern complex is isomorphic to the mapping cone of $\Delta$
shifted by $-1$. Since
$(\iota_p+\overline{\iota_q})\circ\Delta=0$, the zero morphism
$\mathbb C\to 0$, together with
$\iota_p+\overline{\iota_q}$ induces a morphism of complexes from the
Bott--Chern complex to $\mathcal L^\bullet_{p,q}[-1]$.
A direct verification shows that this morphism is a quasi-isomorphism.
\end{remark}
\subsection{Main calculation}

We begin by recalling the notation and explicit bases for the de Rham and Dolbeault cohomology of a two-dimensional theta toroidal group, following \cite{Tan25}. 
Note that the definitions of Bott--Chern and Aeppli cohomology used in \cite{Tan25} differ from the standard ones: in that paper, these groups are defined as quotients of the usual cohomology groups.

For $\tau\in\H=\{z\in\C \mid \Im z >0\}$ and real numbers $p$ and $q$, we define the $\Z$-generated $\R$-independent discrete subgroup $ \Lambda_{\tau,p,q} \subset \C^2$ as follows:  
\begin{align}\label{lattice}
	\Lambda_{\tau,p,q} \coloneqq \left\langle 
	\begin{pmatrix}
	0 \\
	1 \\
	\end{pmatrix},\ 
	\begin{pmatrix}
	1 \\
	p \\
	\end{pmatrix},\ 
	\begin{pmatrix}
	\tau \\
	q \\
	\end{pmatrix}
	\right\rangle.
\end{align}
Let \(X=\mathbb C^2/\Lambda_{\tau,p,q}\) be the resulting complex surface.
The complex surface $X$ is a toroidal group if and only if at least one of \(p\) and \(q\) is irrational \cite{Kop} (see also \cite[Theorem $1.1.4$]{AK}).

Toroidal groups are classified into theta toroidal groups and wild toroidal groups according to number-theoretic conditions, as shown in \cite{Kazama}. 
In particular, \cite[Theorem $4.3$]{Kazama} states that $X=\C^2/\Lambda_{\tau,p,q}$ is a theta toroidal group if and only if there exist $C>0$ and $0<\delta<1$ such that
\begin{align}\label{theta}
\dist (\Z^2, (np,nq)) \geq C\delta^n
\end{align}
for any $n\in\Z_{>0}$, where $\dist$ is the Euclidean distance.

This Diophantine condition ensures that the Dolbeault cohomology groups are finite-dimensional.

The de Rham and Dolbeault cohomologies are described as follows.
\begin{lemma}[\cite{Ume}]\label{dexact_rmk}
Let $X$ be a $2$-dimensional toroidal group. 
For every $d$-closed $C^\infty$ $k$-form $\phi$ on $X$, there exist a unique $k$-form on $X$ with constant coefficients, and a $C^\infty$ $(k-1)$-form $\psi$ on $X$ such that
\begin{align*}
\phi=  \sum \C \bigwedge^k \{dz_1, dz_2, d\bar{z}_1 \} + d\psi.
\end{align*}

\end{lemma}

\begin{lemma}[{\cite[Theorem $4.3$]{Kazama}(see also \cite[ Proposition $2.2.5$]{AK}.)}]\label{dol_repre}
Let $X= \C^2/\Lambda_{\tau,p,q}$ be a $2$-dimensional theta toroidal group. 
Then the Dolbeault cohomology is generated by wedge products of the classes of
\begin{align*}
 \C\bigwedge^p\{dz_1,dz_2\} \wedge \bigwedge^q\{d\bar{z}_1\}.
\end{align*}
\end{lemma}
We now compute the Aeppli cohomology groups.
\begin{theorem}[{Aeppli cohomology, \cite[Theorem 1.3]{Tan25}}]\label{Aeppli}
The Aeppli numbers of a $2$-dimensional non-compact theta toroidal group $X$ are as follows:
\begin{align*}
\begin{array}{ccccc}
     &  & h^{0,0}_{A}=2 &  &  \\
     & h^{1,0}_{A}=2 &  & h^{0,1}_{A}=2 &  \\
    h^{2,0}_{A}=0 &  & h^{1,1}_A=2 &  & h^{0,2}_{A}=0 \\
     & h^{2,1}_{A}=0 &  & h^{1,2}_{A}=0 &  \\
     &  & \phantom{\;.} h^{2,2}_{A}=0 \;. &  & 
\end{array}
\end{align*}
\end{theorem}
The theorem follows by computing the groups \(H_A^{p,q}(X)\) individually in the following lemmas. Throughout this section, \(X\) denotes a two-dimensional theta toroidal group.
\begin{lemma}We have
$H^{0,0}_{A}(X) \simeq \C^2$.
\end{lemma}
\begin{proof}By definition
$$H^{0,0}_{A}(X)=\H^{1}(\cB_{1,1}^\bullet)$$
and the associated long exact sequence is
$$H^{0}(\C)=\C \to H^0(\cO_X \oplus \overline{\cO_X})=\C^2 \to H^{0,0}_{A}(X) \to H^1(\C)=\C^3 \to H^1(\cO_X \oplus \overline{\cO_X})=\C^2. $$
The first map $H^{0}(\C)=\C \to H^0(\cO_X \oplus \overline{\cO_X})=\C^2$ sends the constant function \(1\) to \((1,-1)\), so its cokernel is isomorphic to \(\mathbb C\).

The second map $H^1(\C)=\C^3 \to H^1(\cO_X \oplus \overline{\cO_X})=\C^2$
maps the classes of $d\bar{z}_1, dz_1,dz_2$ to the classes of $(d\bar{z}_1,0),(0,-dz_1),(0,-dz_2)$ respectively.
Since the class of $(0,-dz_2)$ vanishes, its kernel is one-dimensional.

Thus $H^{0,0}_{A}(X)$ is isomorphic to $\C^2$.
An explicit basis is given by \(1\) and \(\Im(z_2)\).
\end{proof}
\begin{lemma}
\label{boundary}
The boundary map $H^{1}(\C)\to \H^1(\cO_X\to \Omega^1_X  )$ is an isomorphism. 
\end{lemma}
\begin{proof}
We first recall an explicit description of the boundary map, following \cite[Chap. VI]{these} (see also \cite{Wu23}). 

We use the notation of \cite[Chap. VI]{these}.
Let $K^{\bullet}$ be a complex of sheaves on the space $X$. We denote by $\{F^pK^{\bullet}\}$ the stupid filtration which does not preserve the cohomology at degree $p$: if $q \geq p$, $F^pK^{q}=K^q$, otherwise $F^pK^{q}=0$.
The corresponding quotient complex is denoted by $\sigma_pK^{\bullet}=K^{\bullet} /F^pK^{\bullet}$.

Let $i: \C \to \sigma_p\Omega^{\bullet} \oplus \sigma_{q}\overline{\Omega^{\bullet}}$ be the morphism of complexes (appearing in the Bott--Chern complex) defined by the diagonal map sending $\C$ into $\cO_X \oplus \overline{\cO_X}$ at degree 0 with a minus sign at the second component and viewing $\C$ as a complex centred at degree 0. 
With the above notations, the Bott--Chern complex is the mapping cone of $i$ which we denote by $\rC^{\bullet}(i)[-1]$.

We denote by ${\ec_X}^{p,q}$ the sheaf of smooth forms of type $(p,q)$ on $X$. For each $p$, the Dolbeault complex $({\ec_X}^{p,\bullet}, \dbar)$ is a fine resolution of $\Omega^p_X$. By conjugation, $({\ec_X}^{\bullet,q}, \d)$ is a fine resolution of $\overline{\Omega^q_X}$.  In particular, $\sigma_{p,\bullet} {\ec_X}^{\bullet,\bullet}$ (resp. $\sigma_{\bullet,q} {\ec_X}^{\bullet,\bullet}$) is the Cartan-Eilenberg resolution of $\sigma_p\Omega^{\bullet}_X$ (resp. $\sigma_q \overline{\Omega^{\bullet}_X}$). Taking the total complex of the double complex, we obtain a fine resolution
$\sigma_{< p, \bullet} {\ec_X}^{\bullet,\bullet}$ of $\sigma_p\Omega^{\bullet}_X$, where, by a slight abuse of notation, $\sigma_{< p, \bullet} {\ec_X}^{\bullet,\bullet}$ denotes the sum of the sheaves of smooth forms of bidegree $(k,l)$ with $k < p$. 
Similarly, $\sigma_{ \bullet,< q} {\ec_X}^{\bullet,\bullet}$ resolves $\sigma_q \overline{\Omega^{\bullet}_X}$. 

Throughout the paper, we use $\sigma_{< p,\bullet} {\ec_X}^{\bullet,\bullet}$ (resp. $\sigma_{\bullet,< q} {\ec_X}^{\bullet,\bullet}$) as the fine resolution of $\sigma_p\Omega^{\bullet}_X$ (resp. $\sigma_q \overline{\Omega^{\bullet}_X}$) to compute the hypercohomology of these complexes. In particular, when we say that the hypercohomology of $\sigma_p\Omega^{\bullet}_X$ (resp. $\sigma_q \overline{\Omega^{\bullet}_X}$) is generated by the classes of certain global smooth forms, we mean that these classes are represented by the same forms viewed as global sections of $\sigma_{< p,\bullet} {\ec_X}^{\bullet,\bullet}$ (resp. $\sigma_{\bullet,< q} {\ec_X}^{\bullet,\bullet}$).

Under these identifications, the inclusion
$\C \to \sigma_p\Omega_X^\bullet$
(resp. $\C \to \sigma_q\overline{\Omega}_X^\bullet$)
is represented by the natural projection
\[
\mathrm{Pr}_p\colon
\ec_X^\bullet \to
\sigma_{< p,\bullet}\ec_X^{\bullet,\bullet}
\]
(resp.\
\[
\overline{\mathrm{Pr}}_q\colon
\ec_X^\bullet \to
\sigma_{\bullet,< q}\ec_X^{\bullet,\bullet}).
\]
For $p,q\geq 1$, the Bott--Chern complex is therefore quasi-isomorphic to the mapping cone of the morphism
\[
(\mathrm{Pr}_p,-\overline{\mathrm{Pr}}_q)\colon
\ec_X^\bullet \longrightarrow
\sigma_{< p,\bullet}\ec_X^{\bullet,\bullet}
\oplus
\sigma_{\bullet,< q}\ec_X^{\bullet,\bullet}
\]
shifted by $-1$.
Consequently, for any \(k\) and \(p\),
the boundary map from $H^k(\C)$ to $\H^k(\sigma_p \Omega_X^{\bullet})$ 
is induced by the corresponding natural projection on global sections.

$\H^1(\cO_X \to \Omega^1_X  )$ can be calculated by the long exact sequence
$$H^{0}(\cO_X) \to H^0( \Omega^1_X  ) \to \H^1(\cO_X\to \Omega^1_X  ) \to H^1(\cO_X) \to H^1( \Omega^1_X  ). $$
Note that the map $H^{0}(\cO_X) \to H^0( \Omega^1_X  )$ is induced by $\d$.
Since the only holomorphic functions are constants,
the map vanishes and its cokernel is isomorphic to $\C^2$.
Similarly, the map $H^1(\cO_X) \to H^1( \Omega^1_X  )$
is induced by $\d$. Since $H^1(\cO_X) $ is generated by the class of  $d \bar{z}_1$,
the map vanishes and its kernel is isomorphic to $\C$.
Thus $\H^1(\cO_X \to \Omega^1_X  )$ is isomorphic to $\C^3$.
From the above long exact sequence one deduces that it is generated by the classes of $dz_1, dz_2, d\bar{z}_1$.
Since the classes of the same forms define a basis of $H^1(X,\C)$, the natural projection corresponding to $H^{1}(\C)\to \H^1(\cO_X\to \Omega^1_X  )$  induces an isomorphism of hypercohomologies.
\end{proof}
\begin{remark}\label{rem-quasi-isom2}
There are two natural fine resolutions of the Bott--Chern complex: the mapping cone of
$(\mathrm{Pr}_p,-\overline{\mathrm{Pr}}_q)$, shifted by $-1$, and
$\mathcal{L}^\bullet_{p,q}$.
For \(p,q\geq1\), the following commutative diagram induces a natural quasi-isomorphism between these two complexes:
\[
\begin{tikzcd}
\C \arrow[r] \arrow[d, "i"] &
\ec_X^\bullet \arrow[r] \arrow[d, "{(\mathrm{Pr}_p,-\overline{\mathrm{Pr}}_q)}"] &
0 \arrow[d]\\
\sigma_p\Omega_X^\bullet \oplus \sigma_q\overline{\Omega}_X^\bullet
\arrow[r] &
\sigma_{< p,\bullet}\ec_X^{\bullet,\bullet}
\oplus
\sigma_{\bullet,< q}\ec_X^{\bullet,\bullet}
\arrow[r] &
\mathcal{L}_{p,q}^\bullet.
\end{tikzcd}
\]
Here, the map from
$\sigma_{< p,\bullet}\ec_X^{\bullet,\bullet}
\oplus
\sigma_{\bullet,< q}\ec_X^{\bullet,\bullet}$
to $\mathcal{L}_{p,q}^\bullet$ is induced by the sum of the natural projections from the two summands.

The assumption $p,q\geq 1$ is necessary for the commutativity of the diagram. We will use this diagram only to understand how Aeppli cohomology classes can be represented by forms with non-constant pluriharmonic coefficients, and this assumption is automatically satisfied in that context. Alternatively, one can use hypercocycle representatives, although these may in general have many non-zero components.
\end{remark}
\begin{lemma}
\label{lem-2d-1,0}
We have
$H^{1,0}_{A}(X) \simeq \C^2$.
\end{lemma}
\begin{proof}By definition,
$$H^{1,0}_{A}(X)=\H^{2}(\cB_{2,1}^\bullet)$$
and the associated long exact sequence is
$$H^{1}(\C) \to \H^1(\cO_X\to \Omega^1_X  ) \oplus H^1(\bar{\cO}_X) \to H^{1,0}_{A}(X) \to H^2(\C) \to \H^2(\cO_X\to \Omega^1_X  ) \oplus H^2(\bar{\cO}_X). $$

By Lemma \ref{boundary}, the map $H^{1}(\C) \to \H^1(\cO_X\to \Omega^1_X  )$ is an isomorphism.
The map $H^{1}(\C) \to \H^1(\cO_X\to \Omega^1_X  ) \oplus H^1(\bar{\cO}_X) $ is thus injective and the cokernel of $H^{1}(\C) \to \H^1(\cO_X\to \Omega^1_X  ) \oplus H^1(\bar{\cO}_X) $
is isomorphic to $\C$.

On the other hand, $\H^2(\cO_X \to \Omega^1_X  )$ can be computed from the long exact sequence
$$H^{1}(\cO_X) \to H^1( \Omega^1_X  ) \to \H^2(\cO_X\to \Omega^1_X  ) \to H^2(\cO_X)=0. $$
The first map $H^{1}(\cO_X) \to H^1( \Omega^1_X  )$
is induced by $\d$. Since $H^{1}(\cO_X)$ is generated by the class of $d\bar{z}_1$,
the map vanishes and its cokernel is isomorphic to $\C^2$.
Thus $\H^2(\cO_X\to \Omega^1_X  )$ is isomorphic to $\C^2$.
The preceding long exact sequence then shows that a basis of $\H^2(\cO_X \to \Omega^1_X  )$  is given by the classes of $dz_1 \wedge  d\bar{z}_1, dz_2 \wedge d\bar{z}_1$.
Therefore the map $H^{2}(\C)\to \H^2(\cO_X\to \Omega^1_X  )$ has kernel $\C\{dz_1 \wedge dz_2\}$.

Thus $H^{1,0}_{A}(X)$ is isomorphic to $\C^2$.
We can also describe an explicit basis of \(H_A^{1,0}(X)\). The global forms $dz_1$ and $\Im(z_2)dz_1$ form such a basis. The class of $dz_1$ is the image of the class represented by \(dz_1\) in $H^1(\overline{\cO_X})$.
Recall that the map $H^{1,0}_{A}(X) \to H^2(\C)$ is induced by the projection of the mapping cone (between $\ec^\bullet_X$ and $\sigma_{<2, \bullet} \ec_X^{\bullet,\bullet} \oplus \sigma_{\bullet,<1} \ec_X^{\bullet,\bullet}$) onto its $\ec^\bullet_X$-component.
The form $dz_1 \wedge dz_2$ is the image, under this projection, of the following (global) hypercocycle
\[
(dz_1 \wedge dz_2,0,0,-2\sqrt{-1}\Im(z_2)dz_1)
\in \ec^2 \oplus \ec^{0,1} \oplus \ec^{1,0} \oplus \ec^{1,0}.
\]
Under the quasi-isomorphisms with the complex $\mathcal{L}^\bullet$ constructed in Remarks \ref{rem-quasi-isom} and \ref{rem-quasi-isom2}, the hypercocycle corresponds to the cocycle $-2\sqrt{-1}\Im(z_2)dz_1$.
By the above exact sequence, this form defines a non-trivial class in $H^{1,0}_A(X)$ which is linearly independent of the class of $dz_1$.
\end{proof}

\begin{lemma}We have
$H^{2,0}_{A}(X) =0$.
\end{lemma}
\begin{proof}By definition,
$$H^{2,0}_{A}(X)=\H^{3}(\cB_{3,1}^\bullet)$$
and the relevant part of the associated long exact sequence is
$$ H^2(\bar{\cO}_X)  \to H^{2,0}_{A}(X) \to \H^3(\C \to \cO_X \to \Omega^1_X \to \Omega_X^2  ) . $$

Note that
$\H^3(\C \to \cO_X \to \Omega^1_X \to \Omega_X^2  )$ is by definition $\H^3(\cB_{3,0}^\bullet)=H^{3,0}_{BC}(X)$
which vanishes because there are no smooth forms of type \((3,0)\) on a complex surface.
On the other hand, $H^2(\bar{\cO}_X)=0$.
Therefore  $H^{2,0}_{A}(X)=0$.
\end{proof}
\begin{lemma}
\label{lem-2d-1,1}
We have
$H^{1,1}_{A}(X) \simeq \C^2.$
\end{lemma}
\begin{proof}By definition,
$$H^{1,1}_{A}(X)=\H^{3}(\cB_{2,2}^\bullet)$$
and the associated long exact sequence is
\[
\begin{aligned}
H^{2}(\C) &\to \H^2(\cO_X \to \Omega^1_X)
\oplus \H^2(\bar{\cO}_X \to \bar{\Omega}^1_X) \\
&\to H^{1,1}_{A}(X) \to H^3(\C)
\to \H^3(\cO_X \to \Omega^1_X)
\oplus \H^3(\bar{\cO}_X \to \bar{\Omega}^1_X).
\end{aligned}
\]

The long exact sequence
$$H^{2}(\Omega^1_X)=0 \to \H^3(\cO_X \to \Omega^1_X  ) \to H^3( \cO_X  )=0 $$
shows that $\H^3(\cO_X \to \Omega^1_X  )=0$.

Note that the map $H^{2}(\C) \to \H^2(\cO_X \to \Omega^1_X  )$ sends the classes of $dz_1 \wedge d\bar{z}_1,dz_2 \wedge d\bar{z}_1$ to the classes of $dz_1 \wedge d\bar{z}_1,dz_2 \wedge d\bar{z}_1$.
On the other hand, the map $H^{2}(\C) \to  \H^2(\bar{\cO}_X \to \bar{\Omega}^1_X  )$ sends the classes of $dz_1 \wedge d\bar{z}_1,dz_2 \wedge d\bar{z}_1$ to the classes of $-dz_1 \wedge d\bar{z}_1,0$.
Recall that $\H^2(\bar{\cO}_X  \to \bar{\Omega}^1_X  )$ is generated by the classes of $dz_1 \wedge d\bar{z}_1,dz_1 \wedge d\bar{z}_2$.

Moreover, the class of $dz_1 \wedge d\bar{z}_2$ in $\H^2(\bar{\cO}_X  \to \bar{\Omega}^1_X  )$ is the image of the class of $dz_1 \wedge dz_2$ (which represents the same class as $dz_1 \wedge d\bar{z}_2$) in $H^{2}(\C)$.

With respect to the concatenation of the following
bases of the two summands,
\[
\bigl(dz_1\wedge d\bar z_1,\ dz_2\wedge d\bar z_1\bigr)
\quad\text{and}\quad
\bigl(dz_1\wedge d\bar z_1,\ dz_1\wedge d\bar z_2\bigr),
\]
the images of the above three generators are
\[
\begin{aligned}
[dz_1\wedge dz_2]&\longmapsto
(0,0,\,0,dz_1\wedge d\bar z_2),\\
[dz_1\wedge d\bar z_1]&\longmapsto
(dz_1\wedge d\bar z_1,0,\,-dz_1\wedge d\bar z_1,0),\\
[dz_2\wedge d\bar z_1]&\longmapsto
(0,dz_2\wedge d\bar z_1,\,0,0).
\end{aligned}
\]
These three images are linearly independent, so the cokernel of $H^{2}(\C) \to \H^2(\cO_X  \to \Omega^1_X  ) \oplus \H^2(\bar{\cO}_X  \to \bar{\Omega}^1_X  )$ is
one-dimensional.

Therefore  $H^{1,1}_{A}(X)$ is isomorphic to $\C^2$.
We can also describe an explicit basis of \(H_A^{1,1}(X)\). The global forms $dz_1 \wedge d\bar{z}_1$, $\Im(z_2)dz_1 \wedge d\bar{z}_1$ form such a basis. 

The class of $dz_1 \wedge d\bar{z}_1$ is the image of the class represented by the same form in $\H^2(\cO_X \to \Omega^1_X  )$.
Recall that the map $H^{1,1}_{A}(X) \to H^3(\C)$ is induced by the projection of the mapping cone (between $\ec^\bullet_X$ and $\sigma_{<2, \bullet} \ec_X^{\bullet,\bullet} \oplus \sigma_{\bullet,<2} \ec_X^{\bullet,\bullet}$ )
onto its $\ec^\bullet_X$-component.
The form $dz_1 \wedge d\bar{z}_1 \wedge dz_2$ is the image, under this projection, of the following (global) hypercocycle
$$(dz_1 \wedge d\bar{z}_1\wedge dz_2,0,0,0,2 \sqrt{-1}\Im(z_2)dz_1 \wedge d\bar{z}_1) \in \ec^3 \oplus \ec^{0,2} \oplus \ec^{1,1} \oplus \ec^{2,0} \oplus \ec^{1,1}.$$
Under the quasi-isomorphisms to the complex $\mathcal{L}^\bullet
$ constructed in Remarks \ref{rem-quasi-isom} and \ref{rem-quasi-isom2},
the hypercocycle corresponds to the cocycle $2 \sqrt{-1}\Im(z_2)dz_1 \wedge d\bar{z}_1 $.
By the above exact sequence, this form defines a non-trivial class in $H^{1,1}_A(X)$ which is linearly independent of the class of $dz_1 \wedge d\bar{z}_1$.
\end{proof}
We point out that the above computation of the cokernel is particularly
simple because the map is represented by a matrix between finite-dimensional
vector spaces of small dimensions. In higher dimensions, although the
corresponding cohomology groups still admit explicit bases consisting of
invariant monomial forms, the computation of the cokernel requires a more
careful analysis: one needs to determine the image of each basis element
and understand how these images are related to one another. This is the
viewpoint we will adopt in the general-dimensional calculation below.
\begin{lemma}We have
$H^{2,1}_{A}(X) =0.$
\end{lemma}
\begin{proof}By definition,
$$H^{2,1}_{A}(X)=\H^{4}(\cB_{3,2}^\bullet)$$
and the relevant part of the associated long exact sequence is
$$H^3(\C) \to \H^3(\cO_X \to \Omega^1_X  \to \Omega^2_X) \oplus \H^3(\bar{\cO}_X \to \bar{\Omega}^1_X  ) \to H^{2,1}_{A}(X) \to H^4(\C)=0 . $$

Since $\C$ is quasi-isomorphic to $\cO_X \to \Omega^1_X  \to \Omega^2_X$,
$\H^3(\cO_X \to \Omega^1_X  \to \Omega^2_X )$ is isomorphic to $H^3(X, \C)=\C$.
Thus the map $H^3(\C) \to \H^3(\cO_X \to \Omega^1_X  \to \Omega^2_X) \oplus \H^3(\bar{\cO}_X \to \bar{\Omega}^1_X  )$ is an isomorphism.
Thus $H^{2,1}_{A}(X)=0$. 
\end{proof}
\begin{lemma}We have
$H^{2,2}_{A}(X) =0.$
\end{lemma}
\begin{proof}By definition,
$$H^{2,2}_{A}(X)=\H^{5}(\cB_{3,3}^\bullet)$$
and the relevant part of the associated long exact sequence is
$$ \H^4(\cO_X \to \Omega^1_X  \to \Omega^2_X) \oplus \H^4(\bar{\cO}_X \to \bar{\Omega}^1_X  \to \bar{\Omega}^2_X) \to H^{2,2}_{A}(X) \to H^5(\C)=0 . $$

Since $\C$ is quasi-isomorphic to $\cO_X \to \Omega^1_X  \to \Omega^2_X$ by the holomorphic Poincaré lemma,
$\H^4(\cO_X \to \Omega^1_X  \to \Omega^2_X )$ is isomorphic to $H^4(X, \C)=0$.
Thus $H^{2,2}_{A}(X)=0$. 
\end{proof}
The Bott--Chern cohomology can be calculated similarly.
\begin{theorem}[{Bott--Chern cohomology, \cite[Theorem 1.1]{Tan25}}]\label{BCcohom}
The Bott--Chern numbers of a $2$-dimensional non-compact theta toroidal group $X$ are as follows:
\begin{align*}
\begin{array}{ccccc}
     &  & h^{0,0}_{BC}=1 &  &  \\
     & h^{1,0}_{BC}=2 &  & h^{0,1}_{BC}=2 &  \\
    h^{2,0}_{BC}=1 &  & h^{1,1}_{BC}=3 &  & h^{0,2}_{BC}=1 \\
     & h^{2,1}_{BC}=1 &  & h^{1,2}_{BC}=1 &  \\
     &  & \phantom{\;.} h^{2,2}_{BC}=0 \;. &  & 
\end{array}
\end{align*}
\end{theorem}
\begin{proof}
The Bott--Chern cohomology groups can be computed in the same way. We record only the relevant long exact sequences.

The long exact sequence
$$H^{0}(\C) \xrightarrow{\simeq} \H^0(\cO_X  ) \to H^{1,0}_{BC}(X) \to H^1(\C) \twoheadrightarrow \H^1(\cO_X )  $$
implies that
$$H^{1,0}_{BC}(X)=\C^2.$$
The long exact sequence
$$H^{1}(\C) \xrightarrow{\simeq} \H^1(\cO_X \to \Omega^1_X ) \to H^{2,0}_{BC}(X) \to H^2(\C) \twoheadrightarrow \H^2(\cO_X \to \Omega^1_X)  $$
implies that
$$H^{2,0}_{BC}(X)=\C.$$
The long exact sequence
$$H^{1}(\C) \twoheadrightarrow H^1(\cO_X \oplus \bar{\cO}_X) \to H^{1,1}_{BC}(X) \to H^2(\C) \to H^2(\cO_X \oplus \bar{\cO}_X)=0  $$
implies that
$$H^{1,1}_{BC}(X)=\C^3.$$
The long exact sequence
$$H^{2}(\C) \twoheadrightarrow  \H^2(\cO_X \to \Omega^1_X ) \oplus H^2(\bar{\cO}_X) \to H^{2,1}_{BC}(X) \to H^3(\C) \to \H^3(\cO_X \to \Omega^1_X)\oplus H^3(\bar{\cO}_X)=0  $$
implies that
$$H^{2,1}_{BC}(X)=\C.$$
The long exact sequence
$$  \H^3(\cO_X \to \Omega^1_X ) \oplus H^3(\bar{\cO}_X \to \bar{\Omega}^1_X)=0 \to H^{2,2}_{BC}(X) \to H^4(\C)=0  $$
implies that
$$H^{2,2}_{BC}(X)=0.$$
\end{proof}
\section{Bott--Chern and Aeppli cohomology of theta toroidal groups}
In this section, we generalize the above calculation to arbitrary-dimensional theta toroidal groups.
\subsection{Preliminaries}
Let $X=\mathbb{C}^n/\Gamma$ be a toroidal group, where $\Gamma\subset\mathbb{C}^n$ is a discrete subgroup of real rank $n+m$ with $0\le m < n$. 
A connected complex Lie group of this form is called toroidal if it admits no nonconstant holomorphic functions. 

By \cite[1.1.12]{AK}, we can choose the following representation of the lattice $\Gamma\subset\mathbb{C}^n$.
In suitable complex coordinates (called toroidal coordinates in \cite{AK}), we may assume that
$\Gamma$ is given by the following matrix 
 $$\begin{pmatrix}
  0 & B\\ 
  \mathrm{Id}_{n-m} & R
\end{pmatrix}$$
where
$B \in \mathrm{Mat}(m,2m; \C)$ and $R \in \mathrm{Mat}(n-m,2m; \R)$.

Note that in these coordinates, although the functions $z_{m+1}, \cdots, z_n$ are not invariant under $\Gamma$,
the imaginary parts $\Im z_{m+1}, \cdots,\Im z_n$ are invariant functions.
In particular,
$dz_{m+1}, \cdots, dz_n$
represent the same de~Rham cohomology classes as $d\bar{z}_{m+1}, \cdots, d\bar{z}_n$.

The de Rham cohomology in toroidal coordinates is described as follows.
\begin{lemma}[\cite{Ume}]\label{dexact_rmk2}
Let $X$ be a toroidal group with the above notation. 
For every $d$-closed $C^\infty$ $k$-form $\phi$ on $X$, there exist a unique $k$-form on $X$ with constant coefficients, and a $C^\infty$ $(k-1)$-form $\psi$ on $X$ such that
\begin{align*}
\phi=  \sum \C \bigwedge^k \{dz_1,\cdots, dz_n,d\bar{z}_1,\cdots, d\bar{z}_m \} + d\psi.
\end{align*}

\end{lemma}

The Dolbeault cohomology is described as follows.
\begin{lemma}[{\cite[Theorem $4.3$]{Kazama}(see also \cite[Proposition $2.2.5$]{AK}.)}]\label{dol_repre2}
Let $X$ be a theta toroidal group as above. 
Then, the Dolbeault cohomology is generated by wedge products of the classes of
\begin{align*}
 \C\bigwedge^p\{dz_1,\cdots, dz_n\} \wedge \bigwedge^q\{d\bar{z}_1, \cdots, d\bar{z}_m\}.
\end{align*}
\end{lemma}
\begin{remark}
A toroidal group that is not a theta toroidal group is called a wild toroidal group. In this case, the Dolbeault cohomology group \(H^{0,1}\) is an infinite-dimensional vector space endowed with a non-Hausdorff topology; see, for example, \cite[Proposition 2.2.7]{AK}. As a consequence, the Hodge decomposition does not hold for a wild toroidal group, and the Hodge--de Rham spectral sequence does not degenerate at the \(E_1\)-page. In particular, the \(\d\dbar\)-lemma fails for wild toroidal groups. An explicit example illustrating this failure is given in \cite[Theorem 3.3 (2)]{KT99}.
\end{remark}
We have the following elementary observation, inspired by
\cite[Remark 3.5]{Wu26}.
\begin{remark}\label{rem-hau-red}
Since the Bott--Chern complex is bounded, the spectral sequence
associated with its hypercohomology degenerates at a finite page.
Moreover, by Lemmas \ref{dexact_rmk2} and \ref{dol_repre2}, all
cohomology groups of the sheaves appearing in the Bott--Chern
complex are finite-dimensional. Hence, all terms on the $E_1$-page
are finite-dimensional, and the same holds for the $E_\infty$-page.
It follows that the hypercohomology groups of the Bott--Chern complex are
finite-dimensional. 
In particular, the Bott--Chern and Aeppli cohomology of the theta toroidal groups are finite-dimensional.

Furthermore, by the arguments in the proofs of Lemmas 1 and 2
of \cite{Ser55}, these cohomology groups are Hausdorff with
respect to their natural topologies, as recalled in
Remark \ref{rem-attention-topo}. In particular, their Hausdorff
reductions coincide with the cohomology groups themselves.
Consequently, Tanaka's computation \cite{Tan25} in fact determines
the full Bott--Chern and Aeppli cohomology groups of
two-dimensional theta toroidal groups.
\end{remark}
As a consequence of Lemma \ref{dol_repre2}, we obtain the following preparation lemma.
This lemma describes the hypercohomology of truncated holomorphic de~Rham complexes in terms of invariant differential forms.
  \begin{lemma} 
\label{lem1}  
Let $X$ be a theta toroidal group as above.  
Let $\mathcal A^\bullet_p$ be the complex \[ \mathcal O \to \Omega^1 \to \cdots \to \Omega^{p-1} \to 0 \] starting in degree $0$. Then $\H^k(X,\mathcal A^\bullet_p)$ can be represented by  invariant constant-coefficient differential forms of bidegree $(r,s)$ with $r \leq p-1, r+s=k$ spanned by $dz_1,\cdots, dz_n,d\bar{z}_1, \cdots, d\bar{z}_m.$
By complex conjugation, $\H^k(X,\overline{\mathcal A}^\bullet_p)$ can be represented by invariant constant-coefficient differential forms of bidegree $(r,s)$ with $s \leq p-1, r+s=k$ spanned by $dz_1,\cdots, dz_m,d\bar{z}_1, \cdots, d\bar{z}_n.$
 \end{lemma} 
\begin{proof} 
Since $\mathcal A^\bullet_1 \simeq \mathcal O_X$, the statement follows for $p=1$ from Lemma \ref{dol_repre2}. 

We now proceed by induction on $p$. Assume that the claim holds for $\mathcal A^\bullet_p$. Consider the short exact sequence of complexes \[ 0 \longrightarrow \Omega^{p}[p] \longrightarrow \mathcal A^{\bullet}_{p+1} \longrightarrow \mathcal A^{\bullet}_{p} \longrightarrow 0, \] where $\Omega^{p}[p]$ denotes the sheaf $\Omega^{p}$ placed in cohomological degree $p$. The associated long exact sequence in hypercohomology yields 
\[\mathbb{H}^{k-1}(X,\mathcal A^{\bullet}_{p}) \to \mathbb{H}^{k-p}(X,\Omega^{p}) \to \mathbb{H}^k(X,\mathcal A^{\bullet}_{p+1}) \to \mathbb{H}^k(X,\mathcal A^{\bullet}_{p}) \to \mathbb{H}^{k-p+1}(X,\Omega^{p}). \]
Note that the boundary maps 
$\mathbb{H}^{k-1}(X,\mathcal A^{\bullet}_{p}) \to \mathbb{H}^{k-p}(X,\Omega^{p})$,
$\mathbb{H}^k(X,\mathcal A^{\bullet}_{p}) \to \mathbb{H}^{k-p+1}(X,\Omega^{p})$
are induced by the $\d$-operator.
Since the representatives are invariant constant-coefficient forms, the $\d$-operator induces the zero map.
In particular,
we have short exact sequences
\[0 \to \mathbb{H}^{k-p}(X,\Omega^{p}) \to \mathbb{H}^k(X,\mathcal A^{\bullet}_{p+1}) \to \mathbb{H}^k(X,\mathcal A^{\bullet}_{p}) \to 0. \]
By the induction hypothesis, $\mathbb{H}^k(X,\mathcal A^{\bullet}_{p})$ can be represented by invariant constant-coefficient differential forms of bidegree $(r,s)$ with $r \leq p-1, r+s=k$ generated by $dz_1,\cdots, dz_n,d\bar{z}_1, \cdots, d\bar{z}_m.$ 
These forms combined with the invariant constant-coefficient differential forms of bidegree $(p,k-p)$ define cohomological classes in $\mathbb{H}^k(X,\mathcal A^{\bullet}_{p+1}) $.
By the short exact sequence, they also give a basis of $\mathbb{H}^k(X,\mathcal A^{\bullet}_{p+1}) $ 
 which completes the induction and hence the proof.   \end{proof} 
 \subsection{Main calculation}
In this subsection, we compute the Bott--Chern and Aeppli cohomology groups of a theta toroidal group.  The computation relies on the long exact sequences in hypercohomology associated with the complexes $\mathcal A_p^\bullet$ and $\mathcal B_{p,q}^\bullet$.  Using Lemma \ref{lem1}, the hypercohomology groups involved can be described in terms of invariant constant-coefficient differential forms.  The dimensions of the Bott--Chern and Aeppli cohomology groups are then obtained by analyzing the kernels and cokernels of the corresponding maps.

Throughout this subsection we represent corresponding cohomology classes by invariant constant-coefficient differential forms of the type
\[ dz_{i_1}\wedge\cdots\wedge dz_{i_a} \wedge d\bar z_{j_1}\wedge\cdots\wedge d\bar z_{j_b} \wedge dz_{k_1}\wedge\cdots\wedge dz_{k_c}, \] where \[ 1\le i_\ell\le m,\qquad 1\le j_\ell\le m,\qquad m+1\le k_\ell\le n, \] and where \(a,b,c\ge 0\). 
The integers \(a,b,c\) denote, respectively, the number of factors \(dz_i\) with \(i\le m\), \(d\bar z_j\) with \(j\le m\), and \(dz_k\) with \(k>m\).
These forms have bidegree \((a+c,b)\) and total degree \(a+b+c\).
Such forms generate the cohomology groups considered below.
 
Let us first consider the Bott--Chern  cohomology group $H^{p,0}_{BC}(X)$ which does not involve the conjugate complex $\mathcal A^{\bullet}_{p}$.
 
\begin{proposition}[Bott--Chern cohomology group of toroidal groups of bidegree $(p,0)$] Let \(X=\mathbb{C}^n/\Gamma\) be a theta toroidal group, where \(\Gamma\subset\mathbb{C}^n\) is a discrete subgroup of rank \(n+m\). Then the Bott--Chern  cohomology group of toroidal groups in bidegree $(p,0)$
satisfies
$$ H^{p,0}_{BC}(X) \simeq  H^{p,0}_{}(X).$$
\end{proposition}
\begin{proof}
Combining the long exact sequences in hypercohomology combining with Lemma \ref{lem1}, we obtain
\[ H^{p-1}(X,\C) \xrightarrow{\simeq}
  \mathbb{H}^{p-1}(X,\mathcal A^{\bullet}_{p}) \to 
 H^{p,0}_{BC}(X) \to
  H^p(X,\C) \twoheadrightarrow
   \mathbb{H}^{p}(X,\mathcal A^{\bullet}_{p}) \]
   which implies that
 $$ H^{p,0}_{BC}(X) \simeq  H^{p,0}_{}(X).$$
\end{proof}
As in Lemma \ref{boundary}, the boundary map between de Rham cohomology and the hypercohomology of $\mathcal{A}_p^\bullet$ is induced by projecting smooth forms onto the appropriate bidegree components.

We first state the main result on the Bott--Chern and Aeppli cohomology of theta toroidal groups; the subsequent discussion justifies these formulas and describes explicit bases.
\begin{theorem}[Bott--Chern and Aeppli Cohomology of theta toroidal groups]
\label{BC-A-theta}
Let \(X=\mathbb{C}^n/\Gamma\) be a theta toroidal group, where \(\Gamma\subset\mathbb{C}^n\) is a discrete subgroup of rank \(n+m\). Then the Bott--Chern and Aeppli cohomology groups of \(X\) are finite-dimensional, and their dimensions are given as follows. 

\medskip \noindent  \emph{(i) Bott--Chern cohomology.} For all \(p,q\ge 1\), \[ \dim_{\mathbb{C}} H^{p,q}_{BC}(X) = \sum_{\substack{a\le p\\ b\le q}} \binom{m}{a}\binom{m}{b} \binom{n-m}{\,p+q-a-b\,}. \] 

\medskip \noindent \emph{(ii) Aeppli cohomology.} For all \(p,q\ge 0\), \[ \dim_{\mathbb{C}} H^{p,q}_{A}(X) = \binom{m}{p}\binom{m}{q} + \sum_{\substack{a\le p\\ b\le q}} \binom{m}{a}\binom{m}{b} \binom{n-m}{\,p+q+1-a-b\,}. \] \end{theorem} 
One can verify the formulas in the preceding example for the case \(n=2,\ m=1\).
All the separate cases have already been discussed in the two-dimensional case.
\begin{proof}
First consider the Aeppli cohomology.
For any $p,q \geq 0$, 
 $$ H^{p,q}_{A}(X) \simeq  \H^{p+q+1}_{}(\cB_{p+1,q+1}^\bullet).$$   
The complex \(\mathcal B_{p+1,q+1}^\bullet\) gives a long exact sequence in hypercohomology
\[ H^{p+q}(\C) \to
  \mathbb{H}^{p+q}(\mathcal A^{\bullet}_{p+1} \oplus\overline{\mathcal A^{\bullet}_{q+1}}) \to 
 H^{p,q}_{A}(X) \to
  H^{p+q+1}(\C) \to
\mathbb{H}^{p+q+1}(\mathcal A^{\bullet}_{p+1} \oplus\overline{\mathcal A^{\bullet}_{q+1}}). \]
As a vector space, \(H^{p,q}_A(X)\) is therefore the direct sum of the cokernel of $H^{p+q}(\C) \to
  \mathbb{H}^{p+q}(\mathcal A^{\bullet}_{p+1} \oplus\overline{\mathcal A^{\bullet}_{q+1}})$
  and
kernel of $  H^{p+q+1}(\C) \to
\mathbb{H}^{p+q+1}(\mathcal A^{\bullet}_{p+1} \oplus\overline{\mathcal A^{\bullet}_{q+1}})$.

We compute these two terms using invariant constant-coefficient differential forms.

First consider the kernel of \[ H^{p+q+1}(\C) \to \mathbb H^{p+q+1}(\mathcal A_{p+1}^\bullet \oplus \overline{\mathcal A_{q+1}^\bullet}). \] This kernel is the intersection of \[ \ker\big(H^{p+q+1}(\C)\to\mathbb H^{p+q+1}(\mathcal A_{p+1}^\bullet)\big) \] and \[ \ker\big(H^{p+q+1}(\C)\to\mathbb H^{p+q+1}(\overline{\mathcal A_{q+1}^\bullet})\big). \]

The first kernel is generated by invariant monomial forms of bidegree \((a+c,b)\) satisfying \[ a+c \ge p+1, \qquad a+b+c=p+q+1. \]
By complex conjugation (of \(dz_{m+1},\dots,dz_n\)), the second kernel is generated by invariant monomial forms of bidegree \((a,b+c)\) satisfying \[ b+c \ge q+1, \qquad a+b+c=p+q+1. \]
Their intersection therefore consists of invariant monomial forms satisfying \[ a\le p, \qquad b\le q, \qquad a+b+c=p+q+1. \] The dimension of this intersection is therefore \[ \sum_{a\le p,\, b\le q} \binom{m}{a}\binom{m}{b}\binom{n-m}{p+q+1-a-b}. \]


Next consider the cokernel of \[ H^{p+q}(\C)\to \mathbb H^{p+q}(\mathcal A_{p+1}^\bullet\oplus\overline{\mathcal A_{q+1}^\bullet}). \] Let \(a,b,c\ge0\) with \(a+b+c=p+q\).
If \(b+c<q\), the corresponding invariant monomial forms represent classes in \(\mathbb H^{p+q}(\overline{\mathcal A_{q+1}^\bullet})\) but not in \(\mathbb H^{p+q}(\mathcal A_{p+1}^\bullet)\). These classes lie in the image of \(H^{p+q}(\C)\), hence contribute trivially to the cokernel. The same conclusion holds for the case \(b+c=q\) and \(a+c>p\). Similarly, if \(a+c<p\), the corresponding classes lie in \(\mathbb H^{p+q}(\mathcal A_{p+1}^\bullet)\) but not in \(\mathbb H^{p+q}(\overline{\mathcal A_{q+1}^\bullet})\), and again lie in the image of \(H^{p+q}(\C)\). The same holds when \(a+c=p\) and \(b+c>q\).
The only remaining case is \[ a+c=p, \qquad b+c=q, \] which implies \(a=p,\ b=q,\ c=0\). These invariant monomial forms represent classes simultaneously in \(\mathbb H^{p+q}(\mathcal A_{p+1}^\bullet)\) and \(\mathbb H^{p+q}(\overline{\mathcal A_{q+1}^\bullet})\). 
The map to \(\mathbb H^{p+q}(\mathcal A_{p+1}^\bullet\oplus\overline{\mathcal A_{q+1}^\bullet})\) is injective on the subspace generated by these
invariant monomial forms. Their images form the diagonal
subspace in the direct sum of the corresponding two
hypercohomology groups. Consequently, the quotient by the
image of \(H^{p+q}(X,\C)\) contributes a vector space of
dimension
\[
\binom{m}{p}\binom{m}{q}.
\]
Consequently \[ h_A^{p,q} = \binom{m}{p}\binom{m}{q} + \sum_{a\le p,\, b\le q} \binom{m}{a}\binom{m}{b}\binom{n-m}{p+q+1-a-b}. \]
As expected, the expression is symmetric in \(p,q\) for Aeppli cohomology.

The following explicit basis forms, representing a basis of the Aeppli cohomology, will be used and will guide the computation of the Aeppli cohomology of wild toroidal groups in the next section.
\begin{remark}
\label{Abase}
A basis of the Aeppli cohomology represented by hypercocycles can be chosen as follows: \[ \Im(z_{k_1})\, dz_{i_1}\wedge\cdots\wedge dz_{i_a} \wedge d\bar z_{j_1}\wedge\cdots\wedge d\bar z_{j_b} \wedge dz_{k_2}\wedge\cdots\wedge dz_{k_c}, \] where \[ 1\le i_1<\cdots<i_a \le m,\qquad 1\le j_1<\cdots<j_b \le m,\qquad m+1\le k_1<\cdots<k_c\le n, \]
with \(a\le p,\ b\le q,\ c=p+q+1-a-b\), together with \[ dz_{i'_1}\wedge\cdots\wedge dz_{i'_p} \wedge d\bar z_{j'_1}\wedge\cdots\wedge d\bar z_{j'_q}, \] where \[ 1\le i'_1<\cdots<i'_p\le m,\qquad 1\le j'_1<\cdots<j'_q\le m. \]
\end{remark}
One can check directly that these forms represent Aeppli cohomology classes. Their linear independence can be established by adapting the arguments in the proofs of Lemmas \ref{lem-2d-1,0} and \ref{lem-2d-1,1}.
More precisely, up to a sign, a form of the first type determines a hypercocycle whose \(\mathcal E_X^{p+q+1}\)-component is its \(\partial\)-differential, one of whose \(\mathcal E_X^{\bullet,\bullet}\)-components is the form itself, and all of whose other components vanish.
The projections onto the \(\mathcal E_X^{p+q+1}\)-component define linearly independent elements of \(H^{p+q+1}(X,\C)\), which generate the kernel of the map
\[
H^{p+q+1}(X,\C)\longrightarrow
\mathbb H^{p+q+1}\bigl(\mathcal A^{\bullet}_{p+1}
\oplus\overline{\mathcal A^{\bullet}_{q+1}}\bigr).
\]
\begin{remark}
\label{Abase2}
A basis of the Aeppli cohomology, represented by $(p,q)$-forms (rather than hypercocycles), can be chosen as follows: \[ \Im(z_{k_1})\, dz_{i_1}\wedge\cdots\wedge dz_{i_a} \wedge d\bar z_{j_1}\wedge\cdots\wedge d\bar z_{j_b} \wedge dz_{k_2}\wedge\cdots\wedge dz_{k_{p-a+1}}
\wedge d\bar z_{k_{p-a+2}}\wedge\cdots\wedge d\bar z_{k_c}, \] where \[  1\le i_1<\cdots<i_a \le m,\qquad 1\le j_1<\cdots<j_b \le m,\qquad m+1\le k_1<\cdots<k_c\le n, \] with \(a\le p,\ b\le q,\ c=p+q+1-a-b\), together with \[ dz_{i'_1}\wedge\cdots\wedge dz_{i'_p} \wedge d\bar z_{j'_1}\wedge\cdots\wedge d\bar z_{j'_q}, \] where \[  1\le i'_1<\cdots<i'_p\le m,\qquad 1\le j'_1<\cdots<j'_q\le m. \]
\end{remark}
\paragraph{} 

The Bott--Chern case is treated similarly. For \(p,q\ge1\), \[ H^{p,q}_{BC}(X)\simeq \mathbb H^{p+q}(\mathcal B_{p,q}^\bullet), \] and the long exact sequence
\[ H^{p+q-1}(\C) \to
  \mathbb{H}^{p+q-1}(\mathcal A^{\bullet}_{p} \oplus\overline{\mathcal A^{\bullet}_{q}}) \to 
 H^{p,q}_{BC}(X) \to
  H^{p+q}(\C) \to
\mathbb{H}^{p+q}(\mathcal A^{\bullet}_{p} \oplus\overline{\mathcal A^{\bullet}_{q}}) \]
shows that \(H^{p,q}_{BC}(X)\) is the direct sum of the corresponding cokernel and kernel. 

The kernel of \[ H^{p+q}(\C)\to \mathbb H^{p+q}(\mathcal A_p^\bullet\oplus\overline{\mathcal A_q^\bullet}) \] consists of invariant monomial forms satisfying \[ a+c\ge p,\qquad b+c\ge q,\qquad a+b+c=p+q. \]
Hence its dimension is \[ \sum_{a\le p,\, b\le q} \binom{m}{a}\binom{m}{b}\binom{n-m}{p+q-a-b}. \] A similar discussion shows that the cokernel of \[ H^{p+q-1}(\C)\to \mathbb H^{p+q-1}(\mathcal A_p^\bullet\oplus\overline{\mathcal A_q^\bullet}) \] vanishes.

Indeed, let \(a,b,c\ge0\) with \(a+b+c=p+q-1\). 
In the Bott--Chern case, the cokernel is generated by invariant monomial forms with either \(b+c \leq q-1\) or \(a+c \leq p-1\); unlike the Aeppli case, no form can simultaneously satisfy \(b+c \leq q-1\) and \(a+c \leq p-1\).
In these cases, the corresponding forms are represented in only one summand, either \(\mathbb{H}^{p+q-1}(\overline{\mathcal{A}_q^\bullet})\) or \(\mathbb{H}^{p+q-1}(\mathcal{A}_p^\bullet)\), and are in the image of \(H^{p+q-1}(\mathbb{C})\). These classes thus contribute trivially to the cokernel. 
Thus \[ h^{p,q}_{BC} = \sum_{a\le p,\, b\le q} \binom{m}{a}\binom{m}{b}\binom{n-m}{p+q-a-b}. \]
As expected, the expression is symmetric in \(p,q\) for Bott--Chern cohomology.
\end{proof}
The following explicit basis forms, representing a basis of the Bott--Chern cohomology, will be used to guide the computation of the Bott--Chern cohomology of wild toroidal groups in the next section. The basis forms in Remark \ref{BCbase} are hypercocycles arising from the hypercohomological description of Bott--Chern cohomology, while those in Remark \ref{BCbase2} represent Bott--Chern cohomology as defined in Definition \ref{BC}.
\begin{remark}
\label{BCbase}
For \(p,q\ge1\), a basis of the Bott--Chern cohomology represented by hypercocycles can be chosen as follows: \[ dz_{i_1}\wedge\cdots\wedge dz_{i_a} \wedge d\bar z_{j_1}\wedge\cdots\wedge d\bar z_{j_b} \wedge dz_{k_1}\wedge\cdots\wedge dz_{k_c}, \] where \[1\le i_1<\cdots<i_a \le m,\qquad 1\le j_1<\cdots<j_b \le m,\qquad m+1\le k_1<\cdots<k_c\le n, \] and \(a\le p,\ b\le q,\ c=p+q-a-b\).
\end{remark}
\begin{remark}
\label{BCbase2}
For $p,q\ge 1$, a basis of the Bott--Chern cohomology, represented by $(p,q)$-forms (rather than hypercocycles), can be chosen as
\[
dz_{i_1}\wedge\cdots\wedge dz_{i_a}
\wedge d\bar z_{j_1}\wedge\cdots\wedge d\bar z_{j_b}
\wedge dz_{k_1}\wedge\cdots\wedge dz_{k_{p-a}}
\wedge d\bar z_{k_{p-a+1}}\wedge\cdots\wedge d\bar z_{k_c},
\]
where
\[
1\le i_1<\cdots<i_a \le m,\qquad 1\le j_1<\cdots<j_b \le m,\qquad m+1\le k_1<\cdots<k_c\le n,
\]
and $a\le p$, $b\le q$, $c=p+q-a-b$.
\end{remark}

\begin{remark}
Using the fact that \(H^{0,1}(X)\) of a wild toroidal group is infinite-dimensional and non-Hausdorff with respect to its natural quotient topology, together with the finite-dimensionality of de Rham cohomology, one can deduce from the long exact sequence that \(H^{1,1}_{BC}(X)\) is also infinite-dimensional.

In general, however, it is unclear whether long exact sequences in hypercohomology, such as the one involving \(\mathbb H^k(\mathcal A_p^\bullet)\), split as topological vector spaces. Moreover, a priori, different choices of fine resolutions of the same complex, each endowed with its natural topology, need not induce isomorphic topological vector spaces upon passing to hypercohomology and equipping it with the quotient topology. In particular, the continuity of a quasi-isomorphism between two such fine resolutions does not necessarily imply that the induced isomorphism on hypercohomology is a homeomorphism.

Consequently, in the following section, we present an alternative method for computing the Bott--Chern and Aeppli cohomology groups of a wild toroidal group, together with their Hausdorff reductions.
Moreover, this provides an alternative, more analytic approach to computing the Bott--Chern and Aeppli cohomology groups of theta toroidal groups.
\end{remark}
\section{Bott--Chern and Aeppli cohomology of wild toroidal groups}
In this section, we compute the Hausdorff reductions of the Bott--Chern and Aeppli cohomology groups of a general toroidal group. Since these groups have already been determined for theta toroidal groups, we focus on the wild case. The same arguments also apply to the theta case and provide an alternative method for computing the corresponding cohomology groups. Compared with the hypercohomological approach, this method is closer in spirit to the original ideas of \cite{Kazama,Ume}.

The organization of this section is as follows. 
In the first subsection, we collect known results on the cohomology groups of wild toroidal groups. In the second subsection, we determine the Hausdorff reductions of the corresponding cohomology groups. In the third subsection, we compute these cohomology groups in certain specific bidegrees without passing to the Hausdorff reduction.
\subsection{Preliminaries}
Since the Bott--Chern and Aeppli cohomology groups of wild toroidal groups are not necessarily Hausdorff, we introduce the notion of Hausdorff reduction (see \cite[Section 2.B]{FRR19}):

\begin{definition}
By a quotient of a Fréchet space (or, briefly, a QF-space), we mean a quotient space $E = F/A$ of a Fréchet space $F$ by an arbitrary linear subspace
$A \subset F$, carrying the induced quotient topology. The Hausdorff reduction of $E$ is then
$$E_{\mathrm{HR}} := E/\overline{\{0\}} = F/\overline{A},$$
 where $\overline{\{0\}}$
  is the smallest closed subspace of $E$, and \(\overline A\) denotes the closure of \(A\) in \(F\).
\end{definition}
When endowed with the uniform topology, 
the Bott--Chern and Aeppli cohomology groups of wild toroidal groups are QF-spaces.
A priori, this topology may differ from that arising from the sheaf-theoretic interpretation of Bott--Chern or Aeppli cohomology endowed with its naturally induced topology, although the resulting topological vector spaces are isomorphic as vector spaces. In the remainder of this section, we do not consider the topology induced by the sheaf-theoretic interpretation.


Let us recall the following computation of Dolbeault cohomology from \cite{FRR19}.
\begin{lemma}[{\cite[Proposition 3.7]{FRR19}}]\label{FRR}
Let $X=\mathbb{C}^n/\Gamma$ be a toroidal group, where $\Gamma\subset\mathbb{C}^n$ is a discrete subgroup of real rank $n+m$ with $0\le m < n$. 
Then, the Hausdorff reduction of the Dolbeault cohomology $H^{p,q}(X)_{\mathrm{HR}}$ admits a basis given by wedge products of the classes of
\begin{align*}
 \C\bigwedge^p\{dz_1,\cdots, dz_n\} \wedge \bigwedge^q\{d\bar{z}_1,\cdots, d\bar{z}_m\}.
\end{align*}
\end{lemma}
For our purposes, we will need the following reformulation of Lemma \ref{FRR}, which can be deduced from the proof of Proposition 3.7 in \cite{FRR19}.

From the general theory of toroidal groups (see e.g.\ \cite[page 51, proof of Proposition 2.2.5]{AK}), any smooth function $f$ on $X$ admits the following Fourier series expansion:
\begin{align}\label{Fourier}
f= \sum_{\sigma\in\Z^{n+m}} a^\sigma (\Im z_{m+1}, \cdots, \Im z_n) \exp\bigl(E^\sigma(\Re z_1, \cdots, \Re z_n, \Im z_{1},\cdots,\Im z_m)\bigr),
\end{align}
where the coefficients $a^\sigma$ are smooth functions, and
\[
E^\sigma(\Re z_1, \cdots, \Re z_n, \Im z_{1},\cdots,\Im z_m)
= \sigma\, E\, (\Re z_1, \cdots, \Re z_n, \Im z_{1},\cdots,\Im z_m)^t
\]
for some complex constant matrix $E$.
\begin{lemma}[{\cite[Proposition 3.7]{FRR19}}]\label{FRR2}
Let $X=\mathbb{C}^n/\Gamma$ be a toroidal group, where $\Gamma\subset\mathbb{C}^n$ is a discrete subgroup of real rank $n+m$ with $0\le m<n$. 
Let $p,q\ge 0$. 

We call an expression of the form
\[
\alpha=\sum_{I,J} f_{I,J}\,dz_I\wedge d\bar z_J,
\]
a formal $(p,q)$-form on $X$
where $I$ and $J$ are multi-indices, and each coefficient $f_{I,J}$ is a formal Fourier series of the form given in \eqref{Fourier}, without any convergence assumption. The operator $\dbar$ is defined termwise on each Fourier component.

Assume that $\alpha$ is $\dbar$-closed. Then there exists a formal $(p,q-1)$-form
\[
\beta=\sum_{I,J'} g_{I,J'}\,dz_I\wedge d\bar z_{J'}
\]
such that 
$$\dbar\beta+\alpha^0=\alpha,$$ 
where each $g_{I,J'}$ is given by a formal Fourier series. Here $\alpha^0$ is uniquely determined by $\alpha$ and is an element of
\begin{align*}
\C\bigwedge^p\{dz_1,\cdots,dz_n\}\wedge \bigwedge^{q}\{d\bar z_1,\cdots,d\bar z_m\} .
\end{align*}
\end{lemma}
By complex conjugation, an analogous result holds for the $\partial$-operator.
All statements concerning formal exactness and the occurrence of monomials are understood mode by mode with respect to the Fourier decomposition.

The following elementary lemma will be used frequently in the next subsection.
\begin{lemma}
\label{exactFourier}
Let $X=\mathbb{C}^n/\Gamma$ be a toroidal group, where $\Gamma\subset\mathbb{C}^n$ is a discrete subgroup of real rank $n+m$ with $0\le m<n$. 
A monomial
\[
dz_{i'_1}\wedge\cdots\wedge dz_{i'_p}
\wedge d\bar z_{j'_1}\wedge\cdots\wedge d\bar z_{j'_q}
\]
with $1\le i'_\ell\le m$ for all $\ell$ does not appear in the
formal Fourier series of any $\d$-exact formal $(p,q)$-form.
Similarly, a monomial of the same form with
$1\le j'_\ell\le m$ for all $\ell$ does not appear in the
formal Fourier series of any $\dbar$-exact formal $(p,q)$-form.
\end{lemma}
\begin{proof}
We argue by contradiction. Suppose that the given monomial appears in the $\sigma=0$ Fourier component of a $\dbar$- or $\partial$-exact $(p,q)$-form. Then there exists a formal $(p,q-1)$-form (respectively $(p-1,q)$-form) whose Fourier expansion satisfies the identity. Taking the $\sigma=0$ Fourier component of this identity and using the fact that $\dbar$ and $\partial$ act termwise, we see that the given monomial appears in the $\sigma=0$ Fourier component of a $\dbar$- or $\partial$-boundary whose coefficients depend only on $\Im z_{m+1},\cdots,\Im z_n$.

For such coefficients, all derivatives with respect to $z_i$ and $\bar z_i$ for $1\le i\le m$ vanish. Hence the operator $\dbar$ (resp. the operator $\partial$) can only produce
terms involving $d\bar z_k$ (resp. $dz_k$) with $k\ge m+1$.

This is a contradiction. Hence the coefficient of the prescribed monomial must vanish.
\end{proof}

\subsection{Main calculation}
We now compute the Hausdorff reduction of the Bott--Chern cohomology of a toroidal group.

We divide the proof into two steps. In Theorem \ref{BCwild}, we construct a family of generators of the corresponding Hausdorff reduction. In Theorem \ref{BCwild2}, we prove their linear independence.
\begin{theorem}
\label{BCwild}
Let $X=\mathbb{C}^n/\Gamma$ be a toroidal group, where $\Gamma\subset\mathbb{C}^n$ is a discrete subgroup of real rank $n+m$ with $0\le m<n$. 
Let $p,q \geq 1$.
Let $\alpha$ be a smooth $(p,q)$-form on $X$ such that $d\alpha=0$. 
Then there exists a formal $(p-1,q-1)$-form $\beta$ such that $\alpha-\partial\bar\partial\beta$ is a finite linear combination of the basis forms described in Remark \ref{BCbase2}. 
\end{theorem}
\begin{proof}
Define $\alpha_{(p,q)}:=\alpha$. Since $\dbar \alpha_{(p,q)}=0$, Lemma \ref{FRR2} yields a formal $(p,q-1)$-form $\gamma_{(p,q-1)}$ such that
\begin{align}
\label{def-alpha}
\alpha_{(p,q)}=\alpha_{(p,q)}^0+\dbar \gamma_{(p,q-1)},
\end{align}
where $\alpha_{(p,q)}^0$ is an invariant constant-coefficient form.
Applying \(\partial\) to this identity, and using \(\partial\alpha_{(p,q)}=0\) together with the fact that \(\alpha_{(p,q)}^0\) has constant coefficients, we obtain \(\partial\bar{\partial}\gamma_{(p,q-1)}=0\).
Applying Lemma \(\ref{FRR2}\) again, we obtain a formal $(p+1,q-2)$-form $\gamma_{(p+1,q-2)}$ such that
\begin{align}
\label{def-gamma}
\d \gamma_{(p,q-1)}=(\d \gamma_{(p,q-1)})^0+\dbar \gamma_{(p+1,q-2)},
\end{align}
where $(\d \gamma_{(p,q-1)})^0$ is an invariant constant-coefficient form.

Proceeding inductively, for each $0\le i\le q-1$ there exist $(p+i,q-i-1)$-forms $\gamma_{(p+i,q-i-1)}$ such that
\[
\d \gamma_{(p+i,q-i-1)}=(\d \gamma_{(p+i,q-i-1)})^0+\dbar \gamma_{(p+i+1,q-i-2)},
\]
where all superscript ${}^0$ terms are invariant constant-coefficient forms, and $\gamma_{(p+q,-1)}=0$.

Note that any monomial in $(\d \gamma_{(p+q-1,0)})^0$ must contain at least one factor $dz_k$ with $k\ge m+1$ by Lemma \ref{exactFourier}. (We will show later, by induction, that such forms appear among the basis elements described in Remark \ref{BCbase}; in particular, they contain at least $q$ factors $dz_k$ with $k\ge m+1$.) 

Therefore, there exists a smooth $(p+q-1,0)$-form $\delta_{(p+q-1,0)}$, whose coefficients are linear functions of $\Im z_{m+1},\cdots,\Im z_n$, such that
\[
\d \delta_{(p+q-1,0)}=(\d \gamma_{(p+q-1,0)})^0,
\]
and $\dbar \delta_{(p+q-1,0)}$ is an invariant constant-coefficient form. 
Note that $\Im z_{m+1},\cdots,\Im z_n$ are globally defined functions on $X$.

More explicitly, up to a sign, $\dbar \delta_{(p+q-1,0)}$ is obtained from $(\d \gamma_{(p+q-1,0)})^0$ by replacing one of the factors among $dz_{m+1},\cdots,dz_n$ with its complex conjugate.

Note that $\gamma_{(p+q-1,0)}-\delta_{(p+q-1,0)}$ is $\d$-closed since $\gamma_{(p+q,-1)}=0$. By Lemma \ref{FRR2}, there exists a formal $(p+q-2,0)$-form $\beta_{(p+q-2,0)}$ such that $\gamma_{(p+q-1,0)}-\delta_{(p+q-1,0)}$ differs from $\d \beta_{(p+q-2,0)}$ by an invariant constant-coefficient form. In particular, taking $\dbar$, we obtain
\[
\dbar \gamma_{(p+q-1,0)}=\dbar \delta_{(p+q-1,0)}+\dbar \d \beta_{(p+q-2,0)}.
\] 
Substituting this into
\[\d \gamma_{(p+q-2,1)}=(\d \gamma_{(p+q-2,1)})^0+\dbar \gamma_{(p+q-1,0)},
\]
we deduce
\begin{align}
\label{start-delta}
\d \gamma_{(p+q-2,1)}=(\d \gamma_{(p+q-2,1)})^0+\dbar \delta_{(p+q-1,0)}+\dbar \d \beta_{(p+q-2,0)}.
\end{align}
We now construct inductively the formal forms $\beta_{(p+q-2-\bullet,\bullet)}$ and the constant-coefficient forms $\delta_{(p+q-1-\bullet,\bullet)}$, with the second index increasing.

(The family of forms $\gamma_{(p+q-1-\bullet,\bullet)}$ is first constructed and then kept fixed. The inductive argument proceeds by canceling its contributions: at each step we introduce correction terms $\delta_{(p+q-1-\bullet,\bullet)}$ and $\beta_{(p+q-2-\bullet,\bullet)}$, which lie in different bidegrees, so as to remove the relevant components of $\gamma_{(p+q-1-\bullet,\bullet)}$. This process is iterated until all such contributions are canceled. The terms $\partial \delta_{(p+q-1-\bullet,\bullet)}$ and $\partial \bar\partial \beta_{(p+q-2-\bullet,\bullet)}$ appear during the induction. 
Only at the final step do they lie in the correct bidegree and contribute to the cohomology and to the Bott--Chern coboundary, respectively.)

Note that any monomial in $(\d \gamma_{(p+q-2,1)})^0+\dbar \delta_{(p+q-1,0)}$ must contain at least one factor $dz_k$ with $k\ge m+1$ by Lemma \ref{exactFourier}. Therefore, there exists a smooth $(p+q-2,1)$-form $\delta_{(p+q-2,1)}$, whose coefficients are linear functions of $\Im z_{m+1},\cdots,\Im z_n$, such that
\begin{align}
\label{def-delta}
\d \delta_{(p+q-2,1)}=(\d \gamma_{(p+q-2,1)})^0+\dbar \delta_{(p+q-1,0)},
\end{align}
and $\dbar \delta_{(p+q-2,1)}$ is an invariant constant-coefficient form.

More explicitly, up to a sign, $\dbar \delta_{(p+q-2,1)}$ is obtained by the following two types of replacements: in a monomial of $(\d \gamma_{(p+q-1,0)})^0$, replace two distinct factors among $dz_{m+1},\cdots,dz_n$ with their complex conjugates; in a monomial of $(\d \gamma_{(p+q-2,1)})^0$, replace one factor among $dz_{m+1},\cdots,dz_n$ with its complex conjugate. In particular, the coefficients of monomials in $(\d \gamma_{(p+q-1,0)})^0$ that do not contain two distinct factors among \(dz_{m+1},\ldots,dz_n\) must vanish.

Note that $\gamma_{(p+q-2,1)}-\delta_{(p+q-2,1)}+\dbar \beta_{(p+q-2,0)}$ is $\d$-closed by Equation (\ref{start-delta}) and (\ref{def-delta}). By Lemma \ref{FRR2}, there exists a formal $(p+q-3,1)$-form $\beta_{(p+q-3,1)}$ such that $\gamma_{(p+q-2,1)}-\delta_{(p+q-2,1)}+\dbar \beta_{(p+q-2,0)}$ differs from $\d \beta_{(p+q-3,1)}$ by an invariant constant-coefficient form. In particular, taking $\dbar$, we obtain
\[
\dbar \gamma_{(p+q-2,1)}=\dbar \delta_{(p+q-2,1)}+\dbar \d \beta_{(p+q-3,1)}.
\]
Substituting this into
\begin{align}
\label{start-gamma}
\d \gamma_{(p+q-3,2)}=(\d \gamma_{(p+q-3,2)})^0+\dbar \gamma_{(p+q-2,1)},
\end{align}
we deduce
\[
\d \gamma_{(p+q-3,2)}=(\d \gamma_{(p+q-3,2)})^0+\dbar \delta_{(p+q-2,1)}+\dbar \d \beta_{(p+q-3,1)}.
\]

Proceeding inductively, we obtain, for each $1\le j\le q-1$, there exist $(p+q-j,j-1)$-forms $\delta_{(p+q-j,j-1)}$ and  $(p+q-j-1,j-1)$-forms $\beta_{(p+q-j-1,j-1)}$ such that
\[
\d \gamma_{(p+q-j-1,j)}=(\d \gamma_{(p+q-j-1,j)})^0+\dbar \delta_{(p+q-j,j-1)}+\dbar \d \beta_{(p+q-j-1,j-1)}.
\]
At the final step of the induction, we use Equation (\ref{def-gamma}) in place of Equation (\ref{start-gamma}).

More explicitly,  up to a sign, the form $\dbar \delta_{(p+q-j,j-1)}$ is obtained as follows. 
For each $j$ and
for each $0\le k\le j$, choose $(j-k)$ distinct factors among $dz_{m+1},\cdots,dz_n$ in a monomial of $(\d \gamma_{(p+q-k-1,k)})^0$ and replace them with their complex conjugates. Each such choice produces a monomial in $\dbar \delta_{(p+q-j,j-1)}$. 
If such a replacement of a monomial is impossible, then the coefficient of this monomial must vanish.

If the replacements are performed in increasing order of the indices, the resulting monomials appear among the basis elements described in Remark \ref{BCbase2}. Moreover, the monomials in $(\d \gamma_{(p+q-k-1,k)})^0$ themselves appear among the basis elements described in Remark \ref{BCbase}.

Note that the above induction procedure does not use the \(\d\)-exactness of \(\d\gamma_{(p+q-j-1,j)}\).
In particular, replacing Equation (\ref{def-gamma}) with Equation (\ref{def-alpha}) at the corresponding step, the same argument shows that there exists a formal $(p-1,q-1)$-form $\beta$ such that $\alpha-\partial\bar\partial\beta$ is a finite linear combination of the basis forms described in Remark \ref{BCbase2}.
\end{proof}
Let us illustrate the argument in the proof of Theorem \ref{BCwild} in the two-dimensional case with the following example.
In this example and in this bidegree, the inductive procedure in the proof of Theorem \ref{BCwild} simplifies and terminates after a single step.
\begin{example}
In this example, we determine generators for the Hausdorff reduction of
$H^{1,1}_{BC}(X)$, where $X$ is a two-dimensional wild toroidal group.

Let $\alpha_{(1,1)}$ be a smooth $d$-closed $(1,1)$-form on $X$.
Then by Lemma \ref{FRR2}, there exist an invariant constant-coefficient form
\[
\alpha_{(1,1)}^0
= C_{1 \bar 1}dz_1 \wedge d \bar z_1
+ C_{2 \bar 1}dz_2 \wedge d \bar z_1
\]
and a formal $(1,0)$-form
\[
\gamma_{(1,0)}
= 2 \sqrt{-1} C_{2 \bar 2}(\Im z_2)\, d z_2
+ 2 \sqrt{-1} C_{1 \bar 2}(\Im z_2)\, d z_1
+ \sum_{\sigma \neq 0}
\bigl(g_{\sigma,1} d z_1 + g_{\sigma,2} d z_2\bigr),
\]
where $C_{2 \bar 2}$ and $C_{1 \bar 2}$ are smooth functions of
$\Im z_2$, and $g_{\sigma,1}$ and $g_{\sigma,2}$ are smooth functions on
$X$, such that
\[
\alpha_{(1,1)}
= \alpha_{(1,1)}^0 + \dbar \gamma_{(1,0)}.
\]
Since $\alpha_{(1,1)}$ is $\d$-closed, $\d \gamma_{(1,0)}$ is
$\dbar$-closed. Then, by Lemma \ref{FRR2}, there exists a constant
$C_{12}$ such that
\[
\d \gamma_{(1,0)}
= C_{12}\, dz_1 \wedge dz_2.
\]
By comparing the terms, we deduce that there exists a constant $K_{12}$
such that
\[
C_{1 \bar 2}(\Im z_2)
= -C_{12}\Im z_2+K_{12},
\]
and
\[
\d \sum_{\sigma \neq 0}
\bigl(g_{\sigma,1} d z_1 + g_{\sigma,2} d z_2\bigr)=0.
\]

Define
\[
\delta_{(1,0)}
= -2 \sqrt{-1} C_{12}\Im z_2\, d z_1
+2 \sqrt{-1} K_{12}\,d z_1.
\]
Then
\[
\d \delta_{(1,0)}
= (\d \gamma_{(1,0)})^0
= \d \gamma_{(1,0)}.
\]
In particular, $\d(\delta_{(1,0)}-\gamma_{(1,0)})=0$.
By the $\d$-analogue of Lemma \ref{FRR2}, there exists a formal function
$\beta_{(0,0)}$ such that
\[
\d \beta_{(0,0)}
= \delta_{(1,0)}-\gamma_{(1,0)}
= -2 \sqrt{-1} C_{2 \bar 2}(\Im z_2)\, d z_2
-\sum_{\sigma \neq 0}
\bigl(g_{\sigma,1} d z_1 + g_{\sigma,2} d z_2\bigr).
\]
In particular, we have
\[
\begin{aligned}
\alpha
={}& C_{1 \bar 1}dz_1 \wedge d \bar z_1
+ C_{2 \bar 1}dz_2 \wedge d \bar z_1
+ \dbar\bigl(
2 \sqrt{-1} C_{1 \bar 2}(\Im z_2)\, d z_1
\bigr)
+ \dbar \d \beta_{(0,0)}\\
={}& C_{1 \bar 1}dz_1 \wedge d \bar z_1
+ C_{2 \bar 1}dz_2 \wedge d \bar z_1
- C_{12}\, d z_1 \wedge d\bar z_2
+ \dbar \d \beta_{(0,0)}.
\end{aligned}
\]
\end{example}

In the zig--zag argument in the proof of Theorem \ref{BCwild}, one repeatedly solves $\dbar$- and $\partial$-equations. As a result, it is not immediate that the forms described in Remark \ref{BCbase2} indeed form a basis of the Hausdorff reduction of the Bott--Chern cohomology. 

To overcome this difficulty, we consider the natural map from Bott--Chern cohomology to de Rham cohomology. The argument is carried out in the following theorem.
\begin{theorem}
\label{BCwild2}
Let $X=\mathbb{C}^n/\Gamma$ be a toroidal group, where $\Gamma\subset\mathbb{C}^n$ is a discrete subgroup of real rank $n+m$ with $0\le m<n$. 
Let $p,q\ge 1$. Then the cohomology classes of the basis forms described in Remark \ref{BCbase2} form a basis of the Hausdorff reduction of $H^{p,q}_{BC}(X)$.
\end{theorem}
For the remaining bidegrees, the Bott--Chern cohomology can be computed explicitly without passing to the Hausdorff reduction. We postpone this computation to Lemma \ref{BCp,0}.
\begin{proof}
By Theorem \ref{BCwild}, the cohomology classes of the basis forms described in Remark \ref{BCbase2} generate the Hausdorff reduction of $H^{p,q}_{BC}(X)$.

The natural map
\[
H^{p,q}_{BC}(X) \longrightarrow H^{p+q}(X,\C)
\]
induced by the identity map on differential forms is continuous.
Since $H^{p+q}(X,\C)$ is Hausdorff, this map factors through a
continuous map
\[
H^{p,q}_{BC}(X)_{\mathrm{HR}}
\longrightarrow H^{p+q}(X,\C).
\]

We now prove that the cohomology classes of the basis forms
described in Remark \ref{BCbase2} are linearly independent.
Consider a linear combination of these forms whose image in
de Rham cohomology vanishes. By Lemma \ref{dexact_rmk2}, its
constant-coefficient representative must then vanish.

On the other hand, since $\Im z_k$ is globally defined for
$k\ge m+1$, the forms $dz_k$ and $d\bar z_k$ represent the same
de Rham cohomology class. Consequently, the forms in
Remark \ref{BCbase2} give distinct monomials in the de Rham
cohomology basis
\[
dz_1,\ldots,dz_n,d\bar z_1,\ldots,d\bar z_m.
\]
Thus all coefficients in the original linear combination
must vanish. This proves the linear independence of the
cohomology classes in the Hausdorff reduction.
\end{proof}
We now compute the Hausdorff reduction of the Aeppli cohomology of a toroidal group using the previous computation of the Bott--Chern cohomology.
\begin{theorem}
\label{Awild}
Let $X=\mathbb{C}^n/\Gamma$ be a toroidal group, where $\Gamma\subset\mathbb{C}^n$ is a discrete subgroup of real rank $n+m$ with $0\le m<n$. 
Let $p,q\ge 0$. Then the cohomology classes of the basis forms described in Remark \ref{Abase2} form a basis of the Hausdorff reduction of $H^{p,q}_{A}(X)$.
\end{theorem}
\begin{proof}
In the following, we assume that at least one of $p,q$ is positive. The case $p=q=0$ can be computed explicitly without passing to the Hausdorff reduction, and we postpone this computation to Lemma \ref{Ae0,0}. By symmetry of Aeppli cohomology, we may assume without loss of generality that $q\ge 1$.

Let $\alpha$ be a $(p,q)$-form such that $\d \dbar \alpha=0$. Since
\[
d(\d \alpha)=\dbar \d \alpha=0,
\]
Theorems \ref{BCwild} and \ref{BCwild2} yield an invariant constant-coefficient $(p+1,q)$-form $\alpha_{\d\text{-}dR}$, as described in Remark \ref{BCbase2}, and a formal $(p,q-1)$-form $\beta$ such that
\[
\d \alpha=\alpha_{\d\text{-}dR}+\d \dbar \beta.
\]
Note that $\alpha_{\d\text{-}dR}$ is uniquely determined by $\d \alpha$.

By Lemma \ref{exactFourier}, every monomial in $\alpha_{\d\text{-}dR}$ contains at least one factor $dz_k$ with $k\ge m+1$. Hence there exists a unique $(p,q)$-form $\alpha_{dR}$, as described in Remark \ref{Abase2}, such that
\[
\d \alpha_{dR}=\alpha_{\d\text{-}dR}.
\]

Explicitly, up to a sign, this amounts to replacing one factor $dz_k$, with $k\ge m+1$, by $\Im z_k$.
If the replacements are made in increasing order of the indices, the resulting
monomials are among the basis elements described in Remark \ref{Abase2}.
(This corresponds, heuristically, to the kernel of the natural map
\[
H^{p,q}_A(X)\longrightarrow H^{p+q+1}(X,\mathbb{C})
\]
appearing in the hypercohomology description in the proof of Theorem \ref{BC-A-theta}.)



Since
\[
\d(\alpha-\alpha_{dR}-\dbar \beta)=0,
\] Lemma \ref{FRR2} yields a a uniquely determined invariant constant-coefficient $(p,q)$-form $\alpha'_{Dol}$ and a formal $(p-1,q)$-form $\gamma$ such that
\[
\alpha-\alpha_{dR}-\dbar \beta=\alpha'_{Dol}+\d \gamma.
\]
Note that $\alpha'_{Dol}$ is uniquely determined by $\alpha-\alpha_{dR}-\dbar \beta$, but not by $\alpha-\alpha_{dR}$ alone.
Moreover, $\alpha'_{Dol}$ is not necessarily a linear combination of the basis forms described in Remark \ref{Abase2}.

The monomials in $\alpha'_{Dol}$ that are not of the type
appearing in Remark \ref{Abase2} are precisely those containing
a factor $d\bar z_k$ with $k\ge m+1$. However, each such monomial can be written as $\dbar$ of the product of the global function $\Im z_k$ with an invariant constant-coefficient form.

Hence these terms can be absorbed into $\dbar \beta$. Consequently,  after possibly replacing $\beta$, there exist an invariant constant-coefficient $(p,q)$-form $\alpha_{Dol}$ and a formal $(p-1,q)$-form $\gamma$ such that
\[
\alpha-\alpha_{dR}-\dbar \beta=\alpha_{Dol}+\d \gamma,
\]
where $\alpha_{dR}+\alpha_{Dol}$ is a linear combination of the basis forms described in Remark \ref{Abase2}.
(Notice that $\alpha_{Dol}$  corresponds, heuristically, to the cokernel of the natural map
\[
\mathbb{H}^{p+q}\bigl(X, \mathcal{A}_{p+1}^\bullet \oplus \overline{\mathcal{A}_{q+1}^\bullet}\bigr)
\longrightarrow H^{p,q}_A(X),
\]
appearing in the hypercohomology description in the proof of Theorem \ref{BC-A-theta}.)

In particular, the cohomology classes of the basis forms described in Remark \ref{Abase2} generate the Hausdorff reduction of the Aeppli cohomology.

We now prove the linear independence of these cohomology
classes. 
We distinguish between the basis forms involving
$\Im z_k$ with $k\ge m+1$ and those that do not.

Consider the continuous map between QF-spaces induced by $\d$
\[
H^{p,q}_{A}(X) \longrightarrow H^{p+q+1}(X,\C),
\]
and the induced continuous map
\[
H^{p,q}_{A}(X)_{HR} \longrightarrow H^{p+q+1}(X,\C).
\]
It follows that the cohomology classes represented by the basis forms involving $\Im z_k$ with $k\ge m+1$ described in Remark \ref{Abase2} are linearly independent in $H^{p,q}_{A}(X)_{HR}$.
Moreover, they are
linearly independent of the classes represented by the
remaining basis forms described in Remark \ref{Abase2}, since the latter lie in the kernel
of the induced map.

Now let $\delta$ be a linear combination of the basis forms
in Remark \ref{Abase2} that do not involve $\Im z_k$ with
$k\ge m+1$, and suppose that its class vanishes in the
Hausdorff reduction of $H^{p,q}_A(X)$.  
This means that
\[
\delta=\d\eta+\dbar\zeta
\]
for some formal $(p-1,q)$-form $\eta$ and formal
$(p,q-1)$-form $\zeta$.

On the other hand, by assumption,  $\delta$ contains no factor $dz_k$ or $d\bar z_k$ with $k\ge m+1$. It follows
from Lemma \ref{exactFourier} that $\delta=0$.  This proves the linear independence in the Hausdorff reduction of the cohomology classes of the basis forms not involving $\Im z_k$ with $k\ge m+1$ described in Remark \ref{Abase2}.
\end{proof}
Let us illustrate the argument in the proof of Theorem \ref{Awild} in the two-dimensional case by the following example.
\begin{example}
In this example, we determine generators for the Hausdorff reduction of $H^{1,1}_{A}(X)$, where $X$ is a two-dimensional wild toroidal group.

Let $\alpha_{(1,1)}$ be a smooth $\d \dbar$-closed $(1,1)$-form on $X$.
Then $\d \alpha_{(1,1)}$ is a smooth $d$-closed $(2,1)$-form.
By Theorem \ref{BCwild2}, there exist a constant $C$ and a formal $(1,0)$-form $\beta$ such that
$$
\d \alpha_{(1,1)} = -\frac{C}{2\sqrt{-1}}\, dz_1 \wedge dz_2 \wedge d\bar z_1 + \d \dbar \beta.
$$
(In the notation of the proof of Theorem \ref{Awild}, we have $\alpha_{\d-dR} = -\frac{C}{2\sqrt{-1}}\, dz_1 \wedge dz_2 \wedge d\bar z_1$.)

It follows that
$$
\d\bigl(\alpha_{(1,1)} - C\, \Im z_2\, dz_1 \wedge d\bar z_1 - \dbar \beta\bigr) = 0.
$$
Hence, by Lemma \ref{FRR2}, there exist a constant $C'$ and a formal $(0,1)$-form $\gamma$ such that
$$
\alpha_{(1,1)} - C\, \Im z_2\, dz_1 \wedge d\bar z_1 - \dbar \beta
= C'\, dz_1 \wedge d\bar z_1 + \d \gamma.
$$
(In the notation of the proof of Theorem \ref{Awild}, $\alpha_{dR} = C\, \Im z_2\, dz_1 \wedge d\bar z_1$ and $\alpha_{Dol}=C'\, dz_1 \wedge d\bar z_1$.)
\end{example}
\begin{remark}
The idea of the proof of Theorem \ref{Awild} is motivated by Proposition \ref{Dem-Sch}, which shows that Aeppli cohomology can be realized as the hypercohomology of the same complex as Bott--Chern cohomology, up to a shift by one in degree.
\end{remark}
\subsection{Complementary results}
In this subsection, we present complementary results that allow us to compute Aeppli and Bott--Chern cohomology explicitly, without passing to the Hausdorff reduction.
\begin{lemma}
\label{Ae0,0}
Let $X=\mathbb{C}^n/\Gamma$ be a toroidal group, where $\Gamma\subset\mathbb{C}^n$ is a discrete subgroup of real rank $n+m$ with $0\le m < n$. 
Then $H^{0,0}_A(X) \simeq \C^{n-m+1}$.
\end{lemma}
\begin{proof}
Consider the following long exact sequence
\[{H}^{0}(X,\C) \to {H}^{0}(X,\cO_X \oplus \overline{\cO}_X) \to {H}^{0,0}_A(X) \to {H}^1(X,\C) \to {H}^{1}(X,\cO_X \oplus \overline{\cO}_X) \]
associated to the short exact sequence
$$0 \to \cO_X \oplus \overline{\cO}_X [1] \to \mathcal B^\bullet_{1,1} \to \C \to 0.$$
By Lemma~\ref{FRR} and Lemma \ref{dexact_rmk2}, the kernel of the natural map \[ H^{1}(X,\C) \longrightarrow H^{1}(X,\cO_X \oplus \overline{\cO}_X) \] has dimension less than or equal to \(n-m\). It follows that \(h^{0,0}_A(X) \le n-m+1\). Since the constant functions and $\Im z_{m+1}, \cdots,\Im z_n$ define $\Gamma$-invariant pluriharmonic functions and  are linearly independent, they form a basis of \(H^{0,0}_A(X)\).
\end{proof}
\begin{lemma}
\label{top-aeppli}
 Let \(X=\mathbb{C}^{n}/\Gamma\) be a toroidal group, where \(\Gamma\subset\mathbb{C}^{n}\) is a discrete subgroup of real rank \(n+m\) with \(0\le m<n\). Then for any \(k\ge 0\), \[ H^{n,k}_{A}(X)=0. \] \end{lemma} 
 \begin{proof}
 For any \(k\ge 0\), the complex \(\mathcal B^{\bullet}_{n+1,k+1}\) used to compute Aeppli cohomology is quasi-isomorphic to the shifted complex \(\mathcal A^{\bullet}_{k+1}[1]\) in the derived category of sheaves of complex vector spaces. 
We therefore have
 \[ H^{n,k}_{A}(X) = \mathbb H^{n+k}(X,\mathcal A^{\bullet}_{k+1}). \] We claim that, for every \(k\ge 0\) and every integer \(l\ge n+k\), \[ \mathbb H^{l}(X,\mathcal A^{\bullet}_{k+1})=0. \] This implies the desired vanishing. We proceed by induction on \(k\). For \(k=0\), the complex \(\mathcal A^{\bullet}_{1}\) is equal to \(\mathcal O_X\) concentrated in degree \(0\). Thus \[ \mathbb H^{l}(X,\mathcal A^{\bullet}_{1}) = H^{l}(X,\mathcal O_X). \] Since \(X\) is a toroidal group of complex dimension \(n\), Malgrange's theorem (see, for example, \cite[Corollary~(4.15), Chapter~IX]{agbook}) implies that \[ H^{l}(X,\mathcal O_X)=0 \qquad \text{for all } l\ge n, \] which establishes the claim for \(k=0\).
Assume now that the claim holds for \(k-1\ge 0\). Consider the short exact sequence of complexes \[ 0 \longrightarrow \Omega^{k}_X[k] \longrightarrow \mathcal A^{\bullet}_{k+1} \longrightarrow \mathcal A^{\bullet}_{k} \longrightarrow 0. \] This induces a long exact sequence in hypercohomology \[ \cdots \longrightarrow H^{l-k}(X,\Omega_X^{k}) \longrightarrow \mathbb H^{l}(X,\mathcal A^{\bullet}_{k+1}) \longrightarrow \mathbb H^{l}(X,\mathcal A^{\bullet}_{k}) \longrightarrow \cdots \] For \(l\ge n+k\),\[ H^{l-k}(X,\Omega_X^{k})=0 \] by Malgrange's theorem. By the induction hypothesis, \[ \mathbb H^{l}(X,\mathcal A^{\bullet}_{k})=0 \qquad \text{for } l\ge n+k. \] Therefore, the exact sequence yields \[ \mathbb H^{l}(X,\mathcal A^{\bullet}_{k+1})=0 \qquad \text{for all } l\ge n+k. \] This completes the induction and proves that \[ H^{n,k}_{A}(X)=0 \qquad \text{for all } k\ge 0. \] 
\end{proof}
\begin{remark}
\label{rem1}
In general, it is difficult to study the hypercohomology groups  \(\mathbb H^{q}(X,\mathcal A^{\bullet}_{p})\) or their Hausdorff reductions for arbitrary \(p\) and \(q\). 

Moreover, the relation between $\mathbb H^{q}(X,\mathcal A^{\bullet}_{p})$ and Dolbeault cohomology is not clear when $X$ is a wild toroidal group.
More precisely,
since the Hausdorff reduction is not an exact functor, there is in general no direct relation between the Hausdorff reduction of \(\mathbb H^{q}(X,\mathcal A^{\bullet}_{p})\) and the Hausdorff reductions of the associated Dolbeault cohomology groups. 
Moreover, Dolbeault cohomology groups may be infinite-dimensional, and a short exact sequence of complexes of vector spaces need not split as a short exact sequence in the category of topological vector spaces. 
 Consequently, even when Dolbeault cohomology is well understood, the structure of \(\mathbb H^{q}(X,\mathcal A^{\bullet}_{p})\) and its Hausdorff reduction remains subtle.
\end{remark}
We now discuss complementary results on Bott--Chern cohomology.
\begin{lemma}
\label{BCp,0}
Let $X=\mathbb{C}^n/\Gamma$ be a toroidal group, where $\Gamma\subset\mathbb{C}^n$ is a discrete subgroup of real rank $n+m$ with $0\le m < n$. 
Then, for any $p$, there is an isomorphism of vector spaces
$$H^{p,0}_{BC}(X) \simeq H^{p,0}(X).$$
\end{lemma}
\begin{proof}
By definition, for any \(p\), the Bott--Chern cohomology group \(H^{p,0}_{BC}(X)\) is a subspace of the Dolbeault cohomology group \(H^{p,0}(X)\). Since \(H^{p,0}(X)\) is represented by invariant constant-coefficient forms, which are both \(\partial\)- and \(\bar\partial\)-closed, they also define classes in \(H^{p,0}_{BC}(X)\). It follows that the natural inclusion induces an isomorphism \[ H^{p,0}_{BC}(X) \simeq H^{p,0}(X). \]
\end{proof}
Note that these cohomology groups are finite-dimensional, and hence the isomorphism is also an isomorphism of topological vector spaces. However, in the following lemma we obtain only an isomorphism of vector spaces, not necessarily of topological vector spaces.
\begin{lemma}
\label{top-bc}
Let $X=\mathbb{C}^n/\Gamma$ be a toroidal group, where $\Gamma\subset\mathbb{C}^n$ is a discrete subgroup of real rank $n+m$ with $0\le m < n$. 
Then, for any \(k\), there is an isomorphism of vector spaces
$$H^{n,k}_{BC}(X) \simeq H^{n,k}(X).$$
\end{lemma}
\begin{proof}
For any \(k\ge 0\), we consider the following exact sequence, as in \cite[Section~3.3]{Var86} \[ \frac{\ker(\partial)\cap\operatorname{Im}(\bar\partial)}{\operatorname{Im}(\partial\bar\partial)} \longrightarrow H^{n,k}_{BC}(X) \longrightarrow H^{n,k}(X) \longrightarrow \frac{\ker(\partial\bar\partial)}{\operatorname{Im}(\bar\partial)+\ker(\partial)} . \] The operator \(\partial\) induces an isomorphism \[ \frac{\ker(\partial\bar\partial)}{\operatorname{Im}(\bar\partial)+\ker(\partial)} \;\xrightarrow{\ \partial\ }\; \frac{\ker(\bar\partial)\cap\operatorname{Im}(\partial)}{\operatorname{Im}(\partial\bar\partial)} . \] In bidegree \((n,k)\), the right-hand side vanishes because there are no forms of bidegree \((n+1,k)\), hence \[ \frac{\ker(\partial\bar\partial)}{\operatorname{Im}(\bar\partial)+\ker(\partial)}=0. \] Similarly, the operator \(\bar\partial\) induces an isomorphism \[ \frac{\ker(\partial\bar\partial)}{\operatorname{Im}(\partial)+\ker(\bar\partial)} \;\xrightarrow{\ \bar\partial\ }\; \frac{\ker(\partial)\cap\operatorname{Im}(\bar\partial)}{\operatorname{Im}(\partial\bar\partial)} . \] On the other hand, by definition there is a natural surjection \[ H^{n,k-1}_{A}(X) = \frac{\ker(\partial\bar\partial)}{\operatorname{Im}(\partial)+\operatorname{Im}(\bar\partial)} \longrightarrow \frac{\ker(\partial\bar\partial)}{\operatorname{Im}(\partial)+\ker(\bar\partial)} . \] 
Since by Lemma \ref{top-aeppli}, \(H^{n,k-1}_{A}(X)=0\), it follows that \[ \frac{\ker(\partial)\cap\operatorname{Im}(\bar\partial)}{\operatorname{Im}(\partial\bar\partial)}=0. \] Therefore, the exact sequence above reduces to an isomorphism \[ H^{n,k}_{BC}(X)\simeq H^{n,k}(X). \]
\end{proof}
A linear isomorphism of vector spaces, even when arising from continuous maps between Fréchet presentations, does not in general induce an isomorphism after Hausdorff reduction.
We do not know whether the above isomorphism is a homeomorphism of topological vector spaces.

As a direct consequence of Malgrange’s theorem, we obtain the following result.
\begin{corollary}
Let $X=\mathbb{C}^n/\Gamma$ be a toroidal group, where $\Gamma\subset\mathbb{C}^n$ is a discrete subgroup of real rank $n+m$ with $0\le m < n$. 
Then 
$$H^{n,n}_{BC}(X) =0.$$
\end{corollary}
  
  Department of Mathematics, Graduate School of Science, Osaka Metropolitan University, 3-3-138 Sugimoto, Osaka 558-8585, Japan.\\
Email address: xiaojun.wu@univ-cotedazur.fr; y25161q@omu.ac.jp.  
  \end{document}